\documentclass[final]{siamltex}

\usepackage{amsrefs,amsmath,amsfonts,amssymb,graphicx,latexsym,color}
\usepackage{hyperref}

\newtheorem{remark}{Remark}

\numberwithin{equation}{section}
\numberwithin{theorem}{section}

\newcommand{\eqdef}{\overset{\mbox{\tiny{def}}}{=}}

\begin{document}

\title{Global Newtonian limit for the Relativistic Boltzmann Equation near Vacuum}


\author{Robert M. Strain\thanks{Princeton University, Department of Mathematics, Fine Hall, Washington Road, Princeton,
NJ 08544 United States (On Leave from the University of Pennsylvania), ({\tt strain at math.upenn.edu}). 
The author has been partly
supported by the NSF fellowship DMS-0602513, and the NSF grant DMS-0901463}}

\maketitle

\begin{abstract}
We study the Cauchy Problem for the relativistic Boltzmann equation with
 near Vacuum initial data. Unique global in time ÒmildÓ solutions are obtained uniformly in the speed of light parameter $c \ge 1$.  We furthermore prove that solutions to the relativistic Boltzmann equation converge to solutions of 
 the Newtonian Boltzmann equation in the limit as $c\to\infty$ on arbitrary time intervals $[0,T]$, with  convergence rate $1/c^{2-\epsilon}$ for any $\epsilon \in(0,2)$. 
This may be the first proof of unique global in time validity of the Newtonian limit for a Kinetic equation.
\end{abstract}



\pagestyle{myheadings}
\thispagestyle{plain}
\markboth{Robert M. Strain}{Global Newtonian limit for the Relativistic Boltzmann Equation}

\section{Introduction}

The dynamics of a relativistic Gas  is modelled by the relativisitic Boltzmann Equation, which can be written as
$$
p^\mu \partial_\mu f = \mathcal{C}(f,f).
$$
The relativistic Boltzmann equation is the central equation in relativistic collisional Kinetic theory.   
Its collision operator is  written in the physics literature \cite{MR1026740,MR635279} as 
\begin{equation}
\mathcal{C}(f,h)
=
 \frac {c}{2}\int_{\mathbb{R}^N} \frac{dq}{q_0}\int_{\mathbb{R}^N}\frac{dq^\prime}{q^\prime_0}\int_{\mathbb{R}^N}\frac{dp^\prime}{p^\prime_0}W(p, q | p^\prime, q^\prime) [f(p^{\prime })h(q^{\prime})-f(p)h(q)].
\label{collisionT}
\end{equation}
The kernel $W(p, q | p^\prime, q^\prime)$ is called the transition rate.  It takes the form
\begin{equation*}
W(p, q | p^\prime, q^\prime) 
=
s\sigma(g, \theta) \delta^{(N+1)}(p^\mu+q^\mu-p^{\mu\prime}-q^{\mu\prime}),
\end{equation*}
where $\sigma(g, \theta)$ is the differential cross-section or scattering kernel; it measure's the interactions between particles.   The rest of the notation is defined just below.

\subsection{History}
First we will discuss a small part of the history of this fundamentally important model.
Lichnerowicz and Marrot \cite{MR0004796} are said to be the first to write down the full relativistic Boltzmann equation, including collisional effects, in 1940.  In 1967, Bichteler \cite{MR0213137} showed that the general relativistic Boltzmann equation has a local solution if the initial distribution function decays exponentially with the energy and if the differential cross-section is bounded.  Dudy{\'n}ski and Ekiel-Je{\.z}ewska \cite{MR933458}, in 1988, showed that the linearized equation admits unique solutions in $L^2$.  Afterwards, Dudy{\'n}ski \cite{MR1031410} studied the long time and small-mean-free-path limits of these solutions.
Then, in 1992, Dudy{\'n}ski and Ekiel-Je{\.z}ewska \cite{MR1151987} proved global existence of large data DiPerna-Lions renormalized solutions \cite{MR1014927}.  In this paper, they use the causality of the relativistic Boltzmann equation, a result which they had previously established 
\cite{MR818441,MR841735}.  This result has since been extended in \cite{MR1714446,MR1676150}.

In 1991, Glassey and Strauss \cite{MR1105532} studied the collision map that carries the pre-collisional momentum of a pair of colliding particles into their momentum post-collision.  
In 1993, Glassey and Strauss \cite{MR1211782} proved existence and uniqueness of smooth solutions which are initially close to a relativistic Maxwellian and in a periodic box.  They also established exponential convergence to Maxwellian.  In 1995, they extended these results to 
the whole space case \cite{MR1321370} where the convergence rate is polynomial.  Analogous results with reduced restrictions 
on the 
cross section were more recently explored in \cite{MR2249574} using the energy method from \cite{MR1908664,MR2000470,MR2013332,guoWS}.  Very recently, the global stability for full range of soft-potentials (in the sense of the assumption in \cite{MR933458})
 and also rapid convergence in a periodic box has been established by the author in \cite{strainSOFT}.

In 1996, Andr{\'e}asson \cite{MR1402446} showed that the gain term is regularizing and the strong $L^1$ compactness.  This is a generalization of Lions \cite{MR1284432,MR1295942} result in the non-relativistic case.  In 1997, Wennberg \cite{MR1480243} proved the regularity of the gain term for both the relativistic and non-relativistic case in a unified framework.   In 2004, Andr{\'e}asson, Calogero and Illner \cite{MR2102321} showed that removal of the loss term for the Boltzmann equation (relativistic or not) leads to finite time blow up of a solution.

In the same year Calogero \cite{MR2098116} proved existence of local-in-time solutions independent of the speed of light and established a rigorous Newtonian limit on that time interval.   Glassey, in 2006, established the first result proving global in time existence of unique mild solutions to the relativistic Boltzmann equation with initial data near the Vacuum state \cite{MR2217287} for certain cross sections.  
Earlier near Vacuum results for the Newtonian Boltzmann equation were already obtained in 1984 by Illner and Shinbrot \cite{MR760333} using in part the iteration method from 1978 by Kaniel and Shinbrot \cite{MR0475532}; for further Newtonian results near Vacuum see also \cite{TartarL,MR1307620,MR776503,MR939503,MR768074,MR996631} and the references therein.
Conditional asymptotic completeness for the relativistic Boltzmann equation was shown in 2007 \cite{MR2459827}; the existence theorem from \cite{MR2217287} does not allow the kind of decay which would be needed for the result of \cite{MR2459827}.
It is these three results 
\cite{MR2098116,MR2217287,MR2459827}
on the relativistic Boltzmann equation that our theorems below can extend.

\subsection{Notation}
Before discussing our main results, we define the problem precisely.
The momentum of a particle is denoted by $p^\mu$, $\mu = 0, 1, \ldots, N$ for $N\ge 2$.  (The most important physical case is of course $N=3$.)
Let the signature of the metric be $(- + \cdots + )$.
We set the rest mass for each particle $m=1$.
The  momentum for each particle is restricted to the mass shell
$
p^\mu p_\mu = -c^2,
$
with
$
p^0>0.
$
The speed of light is a universal constant
$
c  =  299 ~ 792 ~ 458 ~ \text{m/s}
$
in free space.
In this paper we study the Newtonian approximation to Einstein's theory of special relativity, and so we vary the speed of light as a parameter  $c\ge 1$. 
Further with 
$
p \in \mathbb{R}^N,
$
we may write
$
p^\mu=(-p^0, p)
$  
and similarly
$
q^\mu=(-q^0, q).
$  
Thus the energy of a relativistic particle with momentum $p$ is $p_{0}=\sqrt{c^2+|p|^2}$.
The Lorenz inner product is then
$$
p^\mu q_\mu=-p_0 q_0+p\cdot q.
$$
We write the standard inner product for  vectors in $\mathbb{R}^N$ as $p\cdot q =\sum_{i=1}^N p_i q_i$.  

We will now define a few important quantities:
\begin{eqnarray}
s
\eqdef
-(p^\mu+q^\mu) (p_\mu+q_\mu)
=
2\left( -p^\mu q_\mu+c^2\right)\ge 0.
\label{sDEFINITION}
\end{eqnarray}
Above $s$ is the square of the energy in the center-of-momentum system: $p+q=0$.
The relative momentum is denoted
\begin{eqnarray}
g 
\eqdef
\sqrt{(p^\mu-q^\mu) (p_\mu-q_\mu)}
&=&
\sqrt{2(-p^\mu q_\mu-c^2)}.
\label{gDEFINITION}
\end{eqnarray}
Notice that $s=g^2+4c^2$.  This notation is used in \cite{MR635279}, but it is different from some more modern authors by a constant factor.

The angle $\theta$, which is the scattering angle in the center-of-momentum system is defined by
\begin{equation}
\cos\theta= 
(p^\mu - q^\mu) (p_\mu^\prime -q_\mu^\prime)/g^2.
\label{angle}
\end{equation}
For general four-vectors $p^\mu, q^\mu, p^{\mu\prime}, q^{\mu\prime}$ the quantity on the r.h.s. could be unbounded;  however when 
 we assume the collisions are elastic
 $\theta$ can be shown to be well defined \cite[p.113]{MR1379589}.
Momentum conservation for elastic collisions is expressed as
\begin{equation}
p^\mu+q^\mu=p^{\mu\prime}+q^{\mu\prime}.
\label{collisionalCONSERVATION}
\end{equation}
One may also write \eqref{collisionalCONSERVATION}  as 
\begin{eqnarray*}
cp^0+c q^0&=&cp^{0\prime}+ c q^{0\prime}
\\
p+q&=&p^\prime+q^\prime,
\end{eqnarray*}
where the first line represents the principle of conservation of energy and the second line represents the conservation of momentum after a binary collision.
Notice that \eqref{collisionalCONSERVATION} and \eqref{sDEFINITION} together imply that
 $s$ is a collision invariant:
 $
 s(p^\mu,q^\mu)= s(p^{\mu\prime},q^{\mu\prime}).
 $
The Lorentz inner product and $g$ are also:
 $
 p^{\mu } q_\mu
 =
p^{\mu \prime} q_\mu^\prime
 $
 and
  $
 g(p^\mu,q^\mu)= g(p^{\mu\prime},q^{\mu\prime}).
 $

We see that the ``transport term'' is a Lorentz inner product:
$$
p^\mu \partial_\mu 
=
p_0\partial_t + p\cdot \nabla_x.
$$
The usual way mathematician's write the relativistic Boltzmann equation is then
\begin{equation}
\partial _{t}f+\hat{p}\cdot \nabla f=\mathcal{Q}(f,f).
\label{RBF}
\end{equation}
Here in comparison with the notation at the top of this paper, $\mathcal{Q}(f,f) = \mathcal{C}(f,f)/p_0$.  This has become the contemporary standard notation in the mathematics literature.

Above we consider $f=f(t,x,p)$ to be a  function of time $t\in [0,\infty)$, space $x\in \mathbb{R}^N$ and momentum $p\in \mathbb{R}^N$.  
   The normalized velocity of a  particle is denoted 
\begin{equation}
\hat{p}=c\frac{p}{p_0}=\frac{p}{\sqrt{1+|p|^2/c^2}}.
\label{normV}
\end{equation}
Steady states of this model are the well known J\"{u}ttner solution's, also known as the relativistic Maxwellian.
We write the normalized relativistic Maxwellian as
\begin{equation}
J(p)\eqdef \frac{\exp\left(-cp_0\right)}{2(2\pi)^{(N-1)/2} cK_2(c^2)}.
\label{juttner}
\end{equation}
  Here $K_2$ is a bessel function,
see
\cite[p.449]{liboffBK}.
Standard references in relativistic Kinetic theory include 
\cite{MR1898707,MR635279,MR1379589,MR0088362,stewart}.

Before discussing our main results, we will explain two different ways to reduce the collision integrals in \eqref{collisionT}.  
These are analagous to the two well known expressions for the post-collisional velocities in the Newtonian Boltzmann theory.  Each expression that we discuss below converges in the Newtonian limit to the corresponding Boltzmann equation in either the so-called $\omega$ or $\sigma$ notation from \cite{MR1942465}.  The first expression below has become the coordinates that are used in the vast majority of mathematically oriented papers in the field.  The second set of coordinates is not new to the physics literature, but it has  not been widely used mathematically.  These are the coordinates we use to prove Theorem \ref{ue}, they make some of our estimates tractable.

\subsection*{Expression from Glassey-Strauss (1993)}
Glassey and Strauss illustrated in \cite{MR1211782} that a reduction of the collision integrals can be performed, without using Lorentz Transformations, to obtain
\begin{equation}
\mathcal{Q}(f,h)
=
\int_{\mathbb{R}^N\times \mathbb{S}^{N-1}}\mathcal{K}_c(p,q,\omega)[f(p^{\prime })h(q^{\prime})-f(p)h(q)]d\omega dq,
\label{collisionGS}
\end{equation}
where the kernel is
$
\mathcal{K}_c(p,q,\omega)\eqdef
\frac{s\sigma_c (g,\theta)}{p_0q_0}  B(p,q,\omega)
$
with
\begin{equation*}
B(p,q,\omega) = c
\frac{(p_0+q_0)^2p_0q_0\left| \omega\cdot \left(\frac{p}{p_0}-\frac{q}{q_0}\right)\right|}{\left[(p_0+q_0)^2-(\omega\cdot[p+q])^2\right]^2}.
\end{equation*}
In this expression, the post-collisional momentum's are given by
\begin{eqnarray}
p^{\prime } &=&p+a(p,q,\omega )\omega
\nonumber
 \\
q^{\prime } &=&q-a(p,q,\omega )\omega,
\label{postCOLLvelGS}
\end{eqnarray}
where
\begin{equation*}
a(p,q,\omega )=\frac{2(p_{0}+q_{0})p_{0}q_{0}\left\{ \omega \cdot \left( 
\frac{q}{q_{0}}-\frac{p}{p_{0}}\right) \right\} }{(p_{0}+q_{0})^{2}-\left\{
\omega \cdot \left( p+q\right) \right\} ^{2}}.
\end{equation*}
And the energies can be expressed as
$p_0^\prime=p_0+N_0$ and
$q_0^\prime=q_0-N_0$
with 
$$
N_0\eqdef \frac{2 \omega\cdot(p+q)\{p_0(\omega \cdot q)-q_0(\omega \cdot p)\}}{(p_0+q_0)^2-\{\omega\cdot (p+q)\}^2}.
$$
This last expression for the post-collisional energies can be derived from the mass-shell condition, see for instance \cite[p.18]{MR1898707}.
These expressions clearly satisfy the collisional conservations 
\eqref{collisionalCONSERVATION}.  
The angle is then defined by plugging these into \eqref{angle}.
The Jacobian \cite{MR1105532} for the pre-post collisional change of variable $(p, q)\to (p^\prime, q^\prime)$ is found with a non-trivial calculation whose result is
\begin{equation*}
\frac{\partial (p^{\prime },q^{\prime })}{\partial (p,q)}=-
\frac{p_{0}^{\prime }q_{0}^{\prime }}{p_{0}q_{0}}.
\end{equation*}
We do not use these coordinates in the proof of Theorem \ref{ue} in particular because of $N_0$ above.  Simply put there are too many angles, $\omega$, in the expression for $N_0$, which makes it more difficult to utilize these expressions in a few of our proofs below.

Now these variables are analogous to non-relativistic Boltzmann variables, the so-called $\omega$ representation:
\begin{eqnarray}
\notag
\bar{p}^{\prime } &=&p+\omega \cdot \left(q-p\right)\omega
 \\
\bar{q}^{\prime } &=&q-\omega \cdot \left(q-p\right)\omega.
\label{omegaREP}
\end{eqnarray}
With these post-collisional velocities, the non-relativistic Boltzmann equation is
\begin{equation}
\partial _{t}f+p\cdot \nabla f=\mathcal{Q}(f,f),
\label{RBFnewt}
\end{equation}
with the collision operator given by
\begin{equation*}
\mathcal{Q}(f,h)
=
\int_{\mathbb{R}^N\times \mathbb{S}^{N-1}} 
\mathcal{K}_\infty(p,q,\omega) [f(\bar{p}^{\prime })h(\bar{q}^{\prime})-f(p)h(q)]d\omega dq.
\end{equation*}
The kernel above is
\begin{equation}
\mathcal{K}_\infty(p,q,\omega)
=
\left| \omega\cdot \left(p-q\right)\right|\sigma_\infty(|p-q|,\theta).
\label{kernelLIM}
\end{equation}
Indeed when $\sigma_c \to \sigma_\infty $ the formal Newtonian limit, $c\to \infty$, of \eqref{collisionGS} with variables \eqref{postCOLLvelGS} is this standard Newtonian Boltzmann equation.

There is another expression for the collision operator, which we now elaborate.

\subsection*{Center-of-Momentum Collision Operator}
Another method of reduction is  described in the physics literature \cite{MR635279}, which uses Lorentz Transformations in the {\it center-of-momentum system} to reduce the delta functions and obtain 
\begin{equation}
\mathcal{Q}(f,h)
=
\int_{\mathbb{R}^N\times \mathbb{S}^{N-1}} ~ v_c ~ \sigma_c (g,\theta ) ~ [f(p^{\prime })h(q^{\prime})-f(p)h(q)]~ d\omega dq.
\label{collisionCM}
\end{equation}
where $v_c=v_c(p,q)$ is the M\o ller velocity given by
\begin{equation}
v_c=
v_c(p,q)\eqdef 
\frac{c}{2}\sqrt{\left| \frac{p}{p_0}-\frac{q}{q_0}\right|^2-\frac{1}{c^2}\left| \frac{p}{p_0}\times\frac{q}{q_0}\right|^2}
=
\frac{c}{4}\frac{ g\sqrt{s}}{p_0 q_0}.
\label{moller}
\end{equation}
The post-collisional momentum in the expression (\ref{collisionCM}) can be written:
\begin{equation}
\begin{split}
p^\prime&=\frac{p+q}{2}+\frac{g}{2}\left(\omega+(\gamma-1)(p+q)\frac{(p+q)\cdot \omega}{|p+q|^2}\right)
\\
q^\prime&=\frac{p+q}{2}-\frac{g}{2}\left(\omega+(\gamma-1)(p+q)\frac{(p+q)\cdot \omega}{|p+q|^2}\right),
\label{postCOLLvelCMsec2}
\end{split}
\end{equation}
where $\gamma =(p_0+q_0)/\sqrt{s}$.
  The energies are then
\begin{eqnarray*}
p^{0\prime}&=&\frac{p^0+q^0}{2}+\frac{g}{2\sqrt{s}}\omega\cdot (p+q)
\\
q^{0\prime}&=&\frac{p^0+q^0}{2}-\frac{g}{2\sqrt{s}}\omega\cdot (p+q).
\end{eqnarray*}
These clearly satisfy \eqref{collisionalCONSERVATION}.  The angle further satisfies
$
\cos\theta = k\cdot \omega,
$
where $k$ is a unit vector.    In fact, \eqref{postCOLLvelCMsec2} can be defined more generally only up to a (non-unique) Lorentz Transformation, $\Lambda$, satisfying \eqref{lw:lorentzS} as
\begin{equation}
\left(\begin{array}{c} p^{0\prime} \\ p'\end{array}\right) = \frac 12\Lambda^{-1}\left(\begin{array}{c} \sqrt{s} \\ g \omega\end{array}\right),
\quad
\left(\begin{array}{c} q^{0\prime} \\ q'\end{array}\right)
= \frac 12\Lambda^{-1}\left(\begin{array}{c} \sqrt{s} \\ -g \omega\end{array}\right).
\label{postC}
\end{equation}
This will work for any of the Lorentz Transformations in the appendix to this paper;  a complete mathematical derivation of these coordinates can be found in \cite{strainPHD}.\footnote{One sometimes find's in the contemporary mathematics literature the collision operator \eqref{collisionCM} being used in conjunction with the previous coordinates \eqref{postCOLLvelGS}.  When things are done this way the cross section, $\sigma(g,\theta)$, would be implicitly redefined to be different from the cross section originating in the transition rate of the expression for the collision operator \eqref{collisionT}.}

These are the relativistic analogoue of the so-called $\sigma$ representation for the Newtonian Boltzmann collision operator.  In the non-relativistic case, these variables are
\begin{eqnarray*}
\bar{p}^\prime&=&\frac{p+q}{2}+\frac12 |p-q|\omega
\nonumber
\\
\bar{q}^\prime&= &\frac{p+q}{2}-\frac12 |p-q|\omega.
\end{eqnarray*}
And the collision operator interaction is
\begin{equation*}
\mathcal{Q}(f,h)
=
\frac12
\int_{\mathbb{R}^N\times \mathbb{S}^{N-1}}  |p-q| ~ \sigma_\infty(|p-q|,\theta) ~ [f(\bar{p}^{\prime })h(\bar{q}^{\prime})-f(p)h(q)]d\omega dq.
\end{equation*}
Indeed, when $\sigma_c \to \sigma_\infty $ the formal Newtonian limit, $c\to \infty$, of \eqref{collisionCM} with variables \eqref{postCOLLvelCMsec2} is this  system.  
One of the main results of our paper is to show that  this  Newtonian limit can be rigorously established on any time interval.

\section{Statement of the Main Results}
\label{mainres:sec}

We consider the relativistic Boltzmann equation \eqref{RBF} in its mild form
\begin{equation}
f_c^{\#}(t,x,p) = f_{0,c}(x, p) + \int_0^t ~ ds ~ \mathcal{Q}_c^{\#}(f_c, f_c)(s, x, p).
\label{mildRBE}
\end{equation}
We are using standard notation to express a solution along characteristics succinctly as $f_c^{\#}(t,x,p) = f_c(t,x+\hat{p}t,p)$. 
The relativistic collision operator $\mathcal{Q}_c(f_c, f_c)(s, x, p)$ represents either \eqref{collisionGS} or \eqref{collisionCM}, and it will be useful to have the index $c$.

In the Newtonian theory, the mild form is given by 
\begin{equation}
f^{\#}(t,x,p) = f_{0}(x, p) + \int_0^t ~ ds ~ \mathcal{Q}^{\#}(f, f)(s, x, p).
\label{mildRBEnewt}
\end{equation}
Here the characteristics have a slightly different formula: $f^{\#}(t,x,p) = f(t,x+pt,p)$. Additionally, the collision operator 
$\mathcal{Q}(f, f)(s, x, p)$
is here represented as in \eqref{RBFnewt}.  To motivate the following relativistic developments, we will now review the usual solution spaces for a near Vacuum solution to \eqref{mildRBEnewt}.

Solutions to the Newtonian equation are measured in the following norm
$$
\| f^{\#} \|_{\infty}
=
\sup_{t,x,p} 
\frac{ | f^{\#}(t,x,p)|}{\rho_{\infty}(x,p)}.
$$
Here the weight that we use is standard
\begin{equation}
\rho_\infty (x, p) = \exp\left( - \alpha  |x|^2  \right) \mu^{\beta}(p), 
\quad \alpha, \beta > 0.
\label{weightINF}
\end{equation}
The non-relativistic Maxwellian is given by
$$
\mu(v) = (2\pi)^{-N/2} e^{-|v|^2/2}.
$$
Results such as 
\cite{MR760333,TartarL,MR1307620,MR776503,MR939503,MR768074,MR0475532,MR996631}
show that
if $0\le f_{0}(x,p)\in C^0(\mathbb{R}^N_x \times \mathbb{R}^N_p)$ initially and furthermore there
exists a positive number $b_0$ such that if 
\begin{equation}
 \frac{f_{0}(x,p)}{ \rho_\infty(x,p)} \le b_0,
 \label{initialLIMIT}
\end{equation}
 then  a unique non-negative global in time solution $f(t,x,p)$ to the mild form of the Cauchy problem 
 \eqref{mildRBEnewt}
 exists, and satisfies for all time the bound
 $
 \| f^\#\|_\infty
\le b_1
 $
 for a positive constant $b_1$.
 We provide analogs of these statements for the relativistic Boltzmann equation in the following developments.

In order to use the dispersion of solutions along the relativistic characteristics of the transport operator, $x+\hat{p}t$, and thereby control the time integral we introduce the following
 weight, which will  scale in the appropriate way as $c$ becomes large:
\begin{equation}
\rho_c (x, p) = \exp\left( - \alpha  p_0 |x|^2/c  \right) J^{\beta}(p), 
\quad \alpha, \beta > 0.
\label{weightA}
\end{equation}
Above $J^{\beta}(p)$ is the relativistic Maxwellian \eqref{juttner} to the power $\beta$.
This weight function is motivated by the following  invariant quantity:
\begin{multline}
c^3 \frac{t^2}{q_0^\prime} 
+ 
\frac{q_0^\prime}{c} \left|x+t\left( \hat{p} -\hat{q}^\prime \right)\right|^2 
+
c^3 \frac{t^2}{p_0^\prime} 
+ 
\frac{p_0^\prime}{c} \left|x+t\left( \hat{p} - \hat{p}^\prime\right)\right|^2 
\\
=
c^3 \frac{t^2}{p_0} + \frac{p_0}{c} |x|^2 
+
c^3 \frac{t^2}{q_0} + \frac{q_0}{c} \left|x+t\left( \hat{p} - \hat{q}\right)\right|^2. 
\label{weightCALC}
\end{multline}
We found this invariant, or re-found it to discover, during discussions with several experts, that it has existed in the folklore for some time.  
It was furthermore recently recorded in   Jiang \cite{MR2378164}.
This invariance follows from the first order energy and momentum \eqref{collisionalCONSERVATION}.  
It follows  by expanding all the terms in the following way
\begin{gather*}
\frac{p_0}{c} \left|x+t\hat{p}\right|^2 
=
\frac{p_0}{c} |x|^2 - 2t x \cdot p + c p_0 t^2 - \frac{ t^2}{p_0} c^3.
\end{gather*}
We will explain more  towards the end of this section;  the use of this invariant \eqref{weightCALC} requires us to introduce a cut-off \eqref{cutB} in the cross section.  The identity \eqref{weightCALC} forms one of the crucial elements of our global existence proof of Theorem \ref{ue} below.

We will use the following norm to measure the size of solutions to \eqref{mildRBE}:
$$
\| f_c^{\#} \|_{c}
=
\sup_{t,x,p} 
\frac{ | f^{\#}_c(t,x,p)|}{\rho_{c}(x,p)}.
$$
The supremum above is over $t\ge 0$, $x,p \in \mathbb{R}^N$.  Throughout the paper we will use the letter $A$ to denote a generic positive constant, which is independent of the speed of light $c\ge 1$.  These constants, $A$, may change from line to line and they will never depend on important parameters in the problem at hand.  \\

\noindent {\bf Hypothesis on the collision kernel:}   
{\it We assume the collision kernel in \eqref{collisionCM} satisfies the growth/decay estimates
\begin{equation}
0\le \sigma_c (\omega,p,q) 
\le 
\left\{ A_1\left(1+ \left(\frac{g}{1+g}\right)^{\alpha_1}\right)+A_2 ~ g^{-\gamma}\right\} ~ \tilde{\sigma}(\omega) .
\label{hypS}
\end{equation}
Above $A_1$ and $A_2$ are non-negative constants and $\tilde{\sigma}(\omega) \le \sigma_1<\infty$.  We allow $\alpha_1\ge 0$,  $0 \le \gamma < -3$.  
  This includes the hard ball assumption \eqref{hardSPHERE} as a special case.
We also restrict the collision kernel to be supported on  \eqref{cutB}.  
} \\

 See  \cite{DEnotMSI} for a physical discussion of general assumptions, which include our conditions. It would also be possible to prove our theorems with some mild singularities in the angle $\tilde{\sigma}(\omega)$.
We are ready to state our main theorems.  \\

\begin{theorem}[{\bf Uniform Existence}]\label{ue}  Consider the mild form of the relativistic Boltzmann equation \eqref{mildRBE} with cross section of type \eqref{hypS} supported on \eqref{cutB}. 
Choose initial values $0\le f_{0,c}(x,p)\in C^0(\mathbb{R}^N_x \times \mathbb{R}^N_p)$ such that
$$
 \frac{f_{0,c}(x,p)}{ \rho_c(x,p)} \le b.
$$
There is a positive number $b_0$ with the property that if $b\le b_0$ then  a unique non-negative global solution $f_{c}(t,x,p)$ to the mild form of the Cauchy problem exists.

This solution satisfies the global in time estimates 
$$
 \frac{f^{\#}_{c}(t,x,p)}{\rho_c(x,p)}\le  \| f^{\#}_c\|_c  \le b_1.
$$
The norms and of course the solution here are highly dependent on the speed of light, but the constant's $b_0$ and $b_1$ 
are explicit  and do not depend upon $c\ge 1$.  \\
\end{theorem}

 Theorem \ref{ue} provides a proof of global existence for unique mild solutions to the relativistic Boltzmann equation \eqref{mildRBE} for any value of the speed of light $c\ge 1$.  
Next we will study the Newtonian Limit of general solutions to \eqref{mildRBE}.  This next result applies to the solutions in Theorem \ref{ue} but it also applies under much more general conditions.
We  consider the space
$$
\| f_c\|_{L^1_p L^\infty_x}
\eqdef
\int_{\mathbb{R}^N} ~ dp ~ \| f_c(p) \|_{L^\infty(\mathbb{R}^N_x)}. 
$$
It is in this norm that we measure convergence to the Newtonian Boltzmann equation \eqref{RBFnewt} in mild form \eqref{mildRBEnewt}.  
For $h\in\mathbb{R}^N$ define the translation operator $\tau_h^x$ as usual by $\tau_h^x f(x,p) = f(x+h,p)$
and similarly define
$\tau_h^p$ by $\tau_h^p f(x,p) = f(x,p+h)$.
We have\\

\begin{theorem}[{\bf Newtonian Limit}]  
\label{NewtonianLimitThm}
Choose the initial data  $f_0,  \{f_{0,c}\} \in L^1_p L^\infty_x$ for all $c\ge 1$ which initially converge for some $A_1 >0$ as
$$
\| f_{0,c} - f_{0}\|_{L^1_p L^\infty_x}\le A_1/c^{k}, \quad \exists k\in (0,2].
$$
Suppose these initial data $\{f_{0,c}\}$ each lead to global in time unique mild solutions to \eqref{mildRBE}, denoted $f_{c}(t,x,p)$, which satisfy the uniform estimate
$
\| f^{\#}_{c}\|_{c}  \le b_2<\infty.
$
And similarly, suppose $f_{0}$ leads to a global in time mild solution to \eqref{mildRBEnewt} satisfying
$\| f^{\#}\|_{\infty}  \le b_2<\infty$.  Here $b_2$ need not be small.

Consider collision kernels, $\{\sigma_{c}\}$ and $\sigma_\infty$, which satisfy $\left| \sigma_{c} \right| \le \sigma_0<\infty$, and 
$$
\left| \sigma_{c} - \sigma_{\infty} \right| \le A_2 (1+|p|^n+ |q|^m ) / c^{k},
\quad 
\exists
m,n\ge 0, \quad k\in (0,2].
$$
Further suppose $\mathcal{K}_\infty(p,q,\omega)$ from \eqref{kernelLIM} is Lipschitz continuous in $p$, which is a condition on $\sigma_\infty$.
We  rule out  concentration by assuming that for the limit
\begin{equation}
\| \tau_  h^x f_{0} - f_{0}  \|_{L^1_p L^\infty_x}
+
\| \tau_  h^p f_{0}  - f_{0}  \|_{L^1_p L^\infty_x}
\le
A_3  |h|,
\quad 
|h| < 1.
\notag
\end{equation}
Then for any fixed $T>0$ (which is allowed to be large) and  $0
\le t \le T$, the solutions $f_c$ to \eqref{mildRBE}  and $f$ to \eqref{mildRBEnewt}
 corresponding to these initial data converge as:
$$
\| f_{c}(t) - f(t)\|_{L^1_p L^\infty_x}\le A(\delta,T)/c^{k-\delta}, \quad  k\in (0,2], 
\quad \forall \delta \in(0, k).
$$
All the positive constants $\sigma_0,b_2, A_1,A_2,A_3, A(\delta, T)$ are uniform in the speed of light parameter $c$ as needed, but all of these constants are allowed to be large.\\
\end{theorem}

\begin{remark}  
Of course, in Theorem \ref{NewtonianLimitThm}  if initially the $\{f_{0,c}\}$ satisfy the assumptions of Theorem \ref{ue}, then the initial  data convergence condition just stated guarantees that \eqref{initialLIMIT} is satisfied for $f_0$ and the global existence of the Newtonian system \eqref{mildRBEnewt} is known, as long as the constant in \eqref{initialLIMIT} is sufficiently small.

Additionally, the assumption $\left| \sigma_{c} \right| \le \sigma_0<\infty$ is not crucial to Theorem \ref{NewtonianLimitThm}.  In fact singularities such as those in \eqref{hypS} could be allowed by modifications to proofs in this paper.  However we leave these issues for future investigations. 
\end{remark}

\begin{remark} 
Previous works on Newtonian Limits for unique mild or clasical solutions to Kinetic equations, e.g. \cite{MR882384,MR2098116,MR2082817,MR870991,MR840744,MR2136192,MR2023010,MR1277935}, are local-in-time in the following sense.  They generally
 prove that there exists a time interval $[0,T]$ upon which the solution exists independent of $c$. 
Their convergence proofs to the Newtonian approximation in general do continue to hold validity as long as the solutions may exist.
Yet these results do not provide a uniform estimate which establishes global in time existence of a unique solution as the speed of light becomes arbitrarily large.
  We show with the proof of Theorem \ref{ue}   that such estimates are available for the relativistic Boltzmann equation.  Then in Theorem \ref{NewtonianLimitThm} we establish that the solutions form a Cauchy sequence  on any time interval.

There are also Newtonian Limit results in the context of  weak solutions \cite{MR2237676,MR1953298}; in particular \cite{MR2237676} holds in the time periodic case and \cite{MR1953298} establishes weak convergence to the Newtonian approximation globally in time for the Vlasov--Darwin system.   
\end{remark}

\begin{remark}  Previous results on the Newtonian limit for the relativistic Boltzmann equation were obtained by Calogero \cite{MR2098116} in 2004.  In this result, which holds for the periodic case $\mathbb{T}^3$, a local existence theorem is proven on an interval which is independent of the speed of light and then convergence is shown in $L^1_pL^\infty_x$ upon the time interval in which solutions are known to exist.  The proof uses the existence of solutions to the limit system, the Newtonian Boltzmann equation with Hard-Sphere kernel, which is known since \cite{MR0475532,MR760333}.  No convergence rate is given.  
\end{remark}

\begin{remark}  A priori asymptotic completeness results in $L^1$ are shown for the relativistic Boltzmann  equation by Ha, Kim, Lee and Noh in \cite{MR2459827} conditional upon global existence of solutions in a functional space with uniform space-velocity decay.  In the notation of this paper, one of the main results in \cite{MR2459827} shows that given a solution $f_{c}(t,x,p)$ to \eqref{mildRBE}  there exists a unique scattering state $f_{+,c}$ such that
$$
\lim_{t\to\infty}\| f_{c}^{\#}(t) - f_{+,c}(t) \|_{L^1_{x,p}} =0.
$$
Theorem \ref{ue} provides the first existence theorem with the needed uniform decay.  
\end{remark}

\begin{remark}  
In 2006, Glassey \cite{MR2217287} gave the first proof of global existence for mild solutions to the relativistic Boltzmann equation \eqref{mildRBE} near Vaccum.  This  paper appeared 22 years after the  corresponding result \cite{MR760333} for the Newtonian Boltzmann equation first appeared.  Let us point out that a crucial step in \cite{MR760333} is the following second order Newtonian symmetry for the variables \eqref{omegaREP}:
\begin{gather*}
 \left|x+tp\right|^2 
+
\left|x+tq\right|^2 
=
\left|x+t\bar{p}^\prime\right|^2 
+
 \left|x+t\bar{q}^\prime\right|^2,
\end{gather*}
which follows from the Newtonian conservation of energy and momentum:
\begin{gather*}
 \left|p\right|^2 
+
\left|q\right|^2 
=
\left|\bar{p}^\prime\right|^2 
+
 \left|\bar{q}^\prime\right|^2
 \\
  p
+
q
=
\bar{p}^\prime 
+
 \bar{q}^\prime. 
\end{gather*}
It was noted in \cite{MR2217287}  and 10 years earlier in \cite{MR1379589}, that this symmetry fails in the relativistic case; which is true because in particular now the energy and momentum are both first order \eqref{collisionalCONSERVATION}.  Glassey's proof proceeded otherwise by introducing a new weight function
$$
\rho_g(x,p) = e^{-\alpha p_0}(1+|x\times p |^2)^{-(1+\delta)},
\quad
0<\delta<1.
$$
Rather than having a cancellation of the post-collisional momentum via symmetry,  the weight function $\rho_g$ cancel's very nicely the worst effect of the post-collisional momentum because of orthogonality.  The condition $\delta <1$ is used to generate the usual 
reproducing property of the non-linear estimates for the weight function.

The set of cross-sections allowed by this result are given by 
$$
\sigma(p,q,\omega) \le \frac{|\omega \cdot (q \times \hat{p})| \tilde{\sigma}(\omega)}{g(1+g^2)^{\delta + 1/2}},
\quad
0<\delta<1,
$$
and
$$
\int_{S^2} d\omega \frac{\tilde{\sigma}(\omega)}{1+|\omega \cdot z|} \le c |z|^{-1}.
$$ 
This assumption in particular does not include the hard-ball case \eqref{hardSPHERE}.  
\end{remark}

On the other hand, Glassey's result \cite{MR2217287} does not require the assumption which we explain just now.  In the proof of Theorem \ref{ue}, we cut-off part of the angular integration domain $\mathbb{S}^{N-1}$ as  in \eqref{cutB} below.  To describe this we recall the identity \eqref{weightCALC}, which is equivalent to the following
\begin{multline*}
 \frac{p_0}{c} \left|x+t\hat{p}\right|^2 
+
 \frac{q_0}{c} \left|x+t\hat{q}\right|^2 
=
\frac{q_0^\prime}{c} \left|x+t\hat{q}^\prime\right|^2 
+ 
\frac{p_0^\prime}{c} \left|x+t\hat{p}^\prime\right|^2
\\
+
c^3\left( \frac{t^2}{q_0^\prime} 
+
 \frac{t^2}{p_0^\prime} 
-
 \frac{t^2}{p_0}
-
 \frac{t^2}{q_0}\right). 
\end{multline*}
Let us carry out this computation more completely, making use of \eqref{collisionalCONSERVATION},
\begin{multline*}
\frac{p_0}{c} \left|x+t\hat{p}\right|^2 
=
\frac{p_0}{c} |x|^2 - 2t x \cdot p + \frac{p_0}{c}|\hat{p}|^2  t^2
\\
=
\frac{p_0}{c} |x|^2 - 2t x \cdot p + c p_0 t^2 - \frac{ t^2}{p_0} c^3
\\
=
\frac{p_0'}{c} |x|^2 - 2t x \cdot p' + c p_0' t^2 - \frac{ t^2}{p_0} c^3
\\
=
\frac{p_0'}{c} |x|^2 - 2t x \cdot p' + \frac{p_0'}{c}|\hat{p}'|^2 t^2 + \frac{ t^2}{p_0'} c^3 - \frac{ t^2}{p_0} c^3
\\
=
\frac{p_0'}{c} \left|x+t\hat{p}'\right|^2 + \frac{ t^2}{p_0'} c^3 - \frac{ t^2}{p_0} c^3.
\end{multline*}
This shows how to obtain the invariance \eqref{weightCALC}.

Since the terms of the order $t^2$ in the parenthesis of the previous display do not have a sign, they are hard to control pointwise.  Alternatively,
 the weak dispersion from the transport operator is unknown and unlikely to induce the kind of decay which would allow the use of this weight without the cut-off that we are going to introduce now.

For a given constant $B> 0$ and a number $0\le a<1$ and $t > 0$, we define
\begin{gather}
h_c = h(x,p,q,t,c) = \frac{B}{t^2} + a \frac{\alpha q_0 \left|x+t\left( \hat{p} - \hat{q}\right)\right|^2/c }{ t^2} > 0.
\label{cutHC}
\end{gather}
We remark that $h_c$ can generally be quite large  for various values of it's arguments.  
Now we define the cut-off set 
\begin{gather}
\mathcal{B}_c 
=
\left\{\omega : c^3 \left(\frac{1}{p_0}+\frac{1}{q_0}-\frac{1}{p_0'} - \frac{1}{q_0'}\right) \ge - h_c \right\}.
\label{cutB}
\end{gather}
As already noted,  we require that the differential cross sections $\sigma_c$ are all supported on this set \eqref{cutB}.  We also point out that 
$$
h_c \to h_\infty
=
\frac{B}{t^2} + a \frac{\alpha  \left|x+t\left(p - q\right)\right|^2 }{ t^2} > 0, 
\quad 
\text{as}
\quad
c\to \infty.
$$
Furthermore

\begin{lemma}
\label{bA}
$
\lim_{c\to\infty} \mathcal{B}_c    = \mathbb{S}^{N-1}.
$
\end{lemma}

We will prove Lemma \ref{bA} in Section \ref{aeic}.  Lemma \ref{bA} shows that the cut-off does not pose a restriction in the Newtonian Limit.

This assumption may be considered as a relative of Grad's angular cutoff for the relativistic Boltzmann equation.  However, in general, for some values of $p$ and $q$ this assumption will cut-off a larger proportion of $\mathbb{S}^{N-1}$ than the assumption of Grad.  We would like to see a future theory which allows collision kernels such as \eqref{hypS} and does not require this cut-off.  
In our opinion  such a future theory may require either (1) the discovery of a new invariant, or (2)  a substantially different method of proof from the original Illner-Shinbrot \cite{MR760333} approach.
Our main motivation in this work  to illustrate that one can prove a global in time Newtonian limit for this Kinetic equation.

The rest of our paper is organized as follows.  In Section \ref{aeic}, we prove a series of pointwise estimates on algebraic quantities which will be necessary in the following sections.  In Section \ref{geind}, we prove the global existence Theorem \ref{ue} using the coordinates \eqref{postCOLLvelCMsec2}. 
And in Section \ref{gnlimit}, we prove the global Newtonian limit Theorem \ref{NewtonianLimitThm} using instead the Glassey-Strauss coordinates \eqref{postCOLLvelGS}.

Then in Appendix A, we give examples of Lorentz Transformations that lead to the coordinate system \eqref{postCOLLvelCMsec2} which we use in Section \ref{geind} of this paper, as described earlier in \eqref{collisionCM}. 
And finally in  Appendix B, we describe some collisional cross sections, $\sigma(g,\theta)$, that can be found in the physics literature.

\section{Algebraic Estimates independent of $c$}\label{aeic}

In this section, we prove a series of pointwise asymptotic estimates that will be useful in the following sections.  
We first prove an estimate which is stronger than Lemma \ref{bA}.  In particular, this lemma allows us to show that the solutions from Theorem \ref{ue} satisfy the Newtonian Limit Theorem \ref{NewtonianLimitThm}, and in the limit condition \eqref{cutB} disappears.  In fact more is true

\begin{lemma} 
\label{aA}
For a fixed collection $p,q,T,B$,
with $0< t\le T$, there is a constant $c_* = c_*(p,q,T,B)$
such that for $c\ge c_*$ we have
$$
1= {\bf 1}_{\mathcal{B}_c}(\omega), 
$$
where ${\bf 1}_{\mathcal{B}_c}(\omega)$ is the indicator function of the set $\mathcal{B}_c$ defined in \eqref{cutB}.
\end{lemma}

Notice that Lemma \ref{aA} implies Lemma \ref{bA}.  More generally, this estimate implies that the solutions from Theorem \ref{ue} satisfy the results of the Newtonian Limit Theorem \ref{NewtonianLimitThm}.  Let us be more specific.  Consider a cross section $\bar{\sigma}_c$ which satisfies the conditions in Theorem \ref{NewtonianLimitThm}.   Then define
$
\sigma_c = \bar{\sigma}_c {\bf 1}_{\mathcal{B}_c}.
$
Now we may construct global solutions with this cross section, $\sigma_c$, using  Theorem \ref{ue}.  Further the final lower bound in the proof of Lemma \ref{aA} will imply that Lemma \ref{nonLINasy} below remains valid for $\sigma_c$ and this is enough to prove the convergence property in Theorem \ref{NewtonianLimitThm}.

\begin{proof}
We start with the condition on the complement of $\mathcal{B}_c$ that
$$
\frac{1}{p_0}+\frac{1}{q_0}-\frac{1}{p_0'} - \frac{1}{q_0'} \le -\frac{h}{c^3}.
$$
Multiply by $p_0 q_0 p_0' q_0'$ to obtain
$$
(q_0  
+
p_0 ) p_0' q_0'
-
p_0 q_0 ( q_0'
+
p_0' )
\le 
-\frac{h}{c^3} p_0 q_0 p_0' q_0'.
$$
Thus
$$
p_0' q_0'
-
p_0 q_0 
\le 
-\frac{h}{c^3} \frac{p_0 q_0 p_0' q_0'}{(p_0  + q_0 ) }.
$$
We use the energies of the coordinates \eqref{postCOLLvelCMsec2}
\begin{eqnarray*}
p^\prime_0&=&\frac{p_0+q_0}{2}+a_c, 
\quad
a_c = \frac{g}{2\sqrt{s}}\omega\cdot (p+q)
\\
q^\prime_0&=&\frac{p_0+q_0}{2}-a_c.
\end{eqnarray*}
Notice
$$
p_0' q_0' 
=  
\left( \frac{p_0+q_0}{2} \right)^2 - a_c^2
= 
\frac{p_0^2}{4}  
+
\frac{q_0^2}{4}  
+
\frac{p_0 q_0}{2}  
-
 a_c^2.
$$
We plug this in to obtain
$$
\frac{p_0^2}{4}  
+
\frac{q_0^2}{4}  
-
\frac{p_0 q_0}{2}  
-
 a_c^2
\le 
-\frac{h}{c^3} \frac{p_0 q_0 p_0' q_0'}{(p_0  + q_0 ) }.
$$
This implies 
$$
-
 a_c^2
\le 
-\frac{h}{c^3} \frac{p_0 q_0 p_0' q_0'}{(p_0  + q_0 ) }
-
\frac{(p_0-q_0)^2}{4}.
$$
Let 
$$
|p+q|\cos\theta = \omega \cdot (p+q).
$$
Equivalently
$$
\cos^2 \theta
\ge 
\frac{4}{|p+q|^2} \frac{h}{c^3} \frac{p_0 q_0 p_0' q_0'}{(p_0  + q_0 ) }  \frac{s}{g^2}
+
\frac{(p_0-q_0)^2}{|p+q|^2}\frac{s}{g^2} \eqdef \ell^2_c.
$$
We now see that $\ell^2_c$ goes to infinity with rate $O(c^2)$,
since $s=g^2 + 4c^2$.  

We will derive a sequence of (sometimes crude) lower bounds for $\ell^2_c$ which allow us to elaborate on this point.  First of all there is still some weak angular dependence in $\ell^2_c$.  However we can bound from below
$$
2p_0' q_0'
=
c^2 \sqrt{1+|p'|^2/c^2}\sqrt{1+|q'|^2/c^2}
\ge 
c p_0.
$$
We will prove this last inequality in Lemma \ref{lowerBB} below.
We therefore have
$$
\ell^2_c 
\ge 
\frac{2}{|p+q|} \frac{h}{c^3} \frac{c p_0^2 q_0 }{(p_0  + q_0 ) }  \frac{s}{g^2}.
$$
For $0<t \le T$, $h\ge \frac{B}{T}$.   Since $s=g^2 + 4c^2$, we may further bound
$$
\ell^2_c 
\ge  
\frac{2}{|p+q|} \frac{B}{T} \frac{ p_0^2 q_0 }{(p_0  + q_0 ) }  \frac{4}{g^2}.
$$
Now 
$$
p_0^2 q_0= c^3 \left(1+\frac{|p|^2}{c^2}\right)\left(1+\frac{|q|^2}{c^2}\right)^{1/2}  \ge c^3
\sqrt{1+\frac{|p|^2}{c^2}} \sqrt{1+\frac{|q|^2}{c^2}},   
$$ 
and 
$$
p_0  + q_0  
\le 
c \left( \sqrt{1+\frac{|p|^2}{c^2}} + \sqrt{1+\frac{|q|^2}{c^2}}    \right)  
\le
2 c \sqrt{1+\frac{|p|^2}{c^2}} \sqrt{1+\frac{|q|^2}{c^2}}.
$$
Thus
$$
\frac{ p_0^2 q_0 }{(p_0  + q_0 ) }  \ge \frac{c^2}{2}. 
$$
Furthermore $g\le |p-q|$.  All these estimates lead to the crude lower bound
$$
\cos^2\theta 
\ge
\ell^2_c 
\ge  
\frac{8}{|p+q| |p-q|} \frac{B}{T} c^2.
$$
This easily establishes the Lemma by showing $\mathcal{B}_c = \mathbb{S}^{N-1}$ for large $c$.
\end{proof}

We remark that this last lower bound is crude for large values of $|p|$ and $|q|$; in particular as can be seen directly from the quantity $\ell^2_c$ above the set $\mathcal{B}_c$ does not rule out large momentum in general.
This Lemma \ref{aA} is our main crucial use of the coordinates \eqref{postCOLLvelCMsec2}.  The rest of the proof of Theorem \ref{ue} is coordinate independent.  

Now we give a series of pointwise asymptotic estimates.

\begin{lemma}
\label{lowerBB}
We have the pointwise estimate
$$
p_0' q_0'
=
c^2 \sqrt{1+|p'|^2/c^2}\sqrt{1+|q'|^2/c^2}
\ge 
c p_0.
$$
And further
$$
|p-q|
\ge 
g
\ge 
c
\frac{ |p-q|}{\sqrt{p_0q_0}}. 
$$
\end{lemma}

\begin{proof}
The  estimates in this lemma are minor refinements of estimates in \cite{MR1211782} (which are done with $c=1$).  Our estimates allow control of the scaling in $c\ge 1$.  

By energy conservation \eqref{collisionalCONSERVATION} we have
$$
\sqrt{1+|p|^2/c^2}
\le
\sqrt{1+|p'|^2/c^2}+\sqrt{1+|q'|^2/c^2}.
$$
Squaring both sides
\begin{multline*}
1+|p|^2/c^2
\le
1+|p'|^2/c^2+1+|q'|^2/c^2
+
2\sqrt{1+|p'|^2/c^2}\sqrt{1+|q'|^2/c^2}
\\
\le
2\left(1+|p'|^2/c^2+1+|q'|^2/c^2\right)
\\
\le
4\left(1+|p'|^2/c^2\right)\left(1+|q'|^2/c^2\right).
\end{multline*}
Taking square roots and multiplying by $c^2$ yields the claimed inequality.

Next, by a difference of squares, we can write \eqref{gDEFINITION} as
\begin{multline*}
g^2= 2\left(\frac{p_0^2q_0^2-(p\cdot q+c^2)^2}{p_0q_0+p\cdot q+c^2} \right) 
=2 \left(\frac{c^2|p|^2+c^2|q|^2+|p|^2|q|^2-(p\cdot q)^2-2c^2 p\cdot q}{p_0q_0+p\cdot q+c^2} \right) 
\nonumber
\\
=2 \left(\frac{c^2 |p-q|^2+|p\times q|^2}{p_0q_0+p\cdot q+c^2} \right).
\end{multline*}
And the estimate follows easily.  The upper bound proof is the same as in \cite{MR1211782}.
\end{proof}

In the next Lemma, we estimate the pointwise asymptotic convergence rates of a few quantities that will need to be controlled in the sequel.

\begin{lemma}
\label{asymptoticBB}
For $c\ge 1$ we have, under the assumptions of Theorem \ref{NewtonianLimitThm}, that
$$
\left| \mathcal{K}_c(p,q,\omega) - \mathcal{K}_\infty(p,q,\omega)\right|
\le 
A\frac{1+|p|^n+|q|^m}{c^k},
\quad
k\in (0,2],
\quad 
\exists m,n \ge 10.
$$
Furthermore
\begin{equation}
\left| \mathcal{K}_\infty(p,q,\omega)\right|
\le 
A\sigma_0
\left(1+|p|+|q|\right).
\label{kernelUPPERbds}
\end{equation}
We have the upper bound for $\mathcal{K}_c$ as
\begin{equation}
\left| \mathcal{K}_c(p,q,\omega)\right|
\le 
A\sigma_0 \left(1+|p| \right)(1+|q|^2)^{5/2}.
\label{kernelUPPERbdsC}
\end{equation}
By symmetry, since $\mathcal{K}_c(p,q,\omega) = \mathcal{K}_c(q,p,\omega)$, we have
\begin{equation}
\left| \mathcal{K}_c(p,q,\omega)\right|
\le 
A\sigma_0 \left(1+|q| \right)(1+|p|^2)^{5/2}.
\notag
\end{equation}
Notice that all of the constants $A>0$ above are independent of $c\ge 1$.
\end{lemma}

\begin{proof}  These basic pointwise estimates were shown for the most part in Calogero \cite{MR2098116}.  In particular 
the first estimate for the difference of $\left| \mathcal{K}_c(p,q,\omega) - \mathcal{K}_\infty(p,q,\omega)\right|$ follows from
\cite[Lemma 1 part (b)]{MR2098116} when combined with the estimate for the difference
$
\left| \sigma_{c} - \sigma_{\infty} \right|
$
 in Theorem \ref{NewtonianLimitThm} and \eqref{kernelUPPERbds}.  Note that $\mathcal{K}_c$ in this paper is defined differently from $\mathcal{K}_c$ in \cite{MR2098116}.
 Furthermore, \eqref{kernelUPPERbds} follows trivially from the assumptions in Theorem \ref{NewtonianLimitThm} (in particular that $\left| \sigma_{\infty} \right| \le \sigma_0<\infty$) and \eqref{kernelLIM}.
 Lastly the estimate \eqref{kernelUPPERbdsC} is proven as an intermediate step in the proof of \cite[Lemma 1 part (c)]{MR2098116}.
\end{proof}

We also have estimates for the difference between the  post-collisional relativistic momentum \eqref{postCOLLvelGS}
and the Newtonian post-collisional momentum \eqref{omegaREP} as
\begin{equation}
\left| \bar{p}^{\prime} - p^{\prime} \right|
+
\left| \bar{q}^{\prime} - q^{\prime} \right| \le A (|p|+|q|)^3/c^2,
\label{postCOLLest}
\end{equation}
This is shown in \cite[Lemma 1 part (a)]{MR2098116}.

Furthermore, for the difference of normalized velocities we have
\begin{gather}
\left| \frac{p}{\sqrt{1+|p|^2/c^2}} - p \right|
\le
A\left(\frac{1+ |p|^3}{c^2} \right).
\label{diffPhat}
\end{gather}
This estimate is proven with a first order taylor expansion of the denominator.  The first terms cancel and the remainder terms are bounded as the upper bound in  \eqref{diffPhat}.




%


\begin{lemma}\label{JSOLbounds} The relativistic Maxwellian  $J(p)$ as in \eqref{juttner} is bounded and has exponential decay independent of the speed of light.  In particular
\begin{equation*}
A^{-1} e^{-\frac12|p|^2} \le J(p)\le A e^{-|p|},
\end{equation*}
where $A>0$ does not depend on the speed of light $c$.
\end{lemma}

\begin{proof}  In (\ref{juttner}) $K_2(x)$ for say $x\ge 1$ is the Bessel function:
\begin{equation*}
K_2(x) 
=\frac{x^2}{3}\int_1^\infty e^{-x s} (s^{2}-1)^{3/2} ds.
\end{equation*}
We translate $s\to s+1$ to obtain
$$
K_2(x)
=\frac{x^2}{3}e^{-x }\int_0^\infty e^{-x s} s^{3/2}(s+2)^{3/2} ds.
$$
And we change variables $s\to s/x$ to obtain
$$
K_2(x)
=\frac{1}{3\sqrt{x}}e^{-x }\int_0^\infty e^{- s} s^{3/2}\left(\frac{s}{x}+2\right)^{3/2} ds.
$$
From the last representation, we conclude that (for $x\ge 1$)
$$
A^{-1} x^{-1/2} e^{-x}\le K_2(x)
\le Ax^{-1/2} e^{-x}.
$$
Plugging this inequality into (\ref{juttner}) we conclude
\begin{equation}
A^{-1} e^{c^2}e^{-cp_0}\le J(p) \le A e^{c^2}e^{-cp_0}.
\label{jutBOUNDS}
\end{equation}
We {\it claim} that 
\begin{eqnarray}
\chi_1(|p|^2/c^2) \eqdef 1+\frac12 \frac{|p|^2}{c^2}-\sqrt{1+|p|^2/c^2}\ge 0
\label{chi1INEQUALTIES}
\\
\chi_2(c) \eqdef c^2-c^2\sqrt{1+|p|^2/c^2}\le \chi_2(1)=1-\sqrt{1+|p|^2}, ~~c\ge 1.
\label{chi2INEQUALTIES}
\end{eqnarray}
Then \eqref{chi1INEQUALTIES} establishes the lower bound and (\ref{chi2INEQUALTIES}) establishes the upper bound. 

We first establish (\ref{chi1INEQUALTIES}), consider $\chi_1(x)$ for $x\ge 0$.   We differentiate
$$
\chi^\prime_1(x)=\frac{1}{2}-\frac{1}{2}(1+x)^{-1/2}\ge 0.
$$
Since $\chi_1(x)$ is increasing and $\chi_1(0)=0$ we have (\ref{chi1INEQUALTIES}).

We finish off the proof by establishing (\ref{chi2INEQUALTIES}).  We will show that $\chi^\prime_2(c)<0$ for any fixed $|p|>0$ and for all $c\ge 0$.  We compute  
$$
\chi^\prime_2(c)=2c-\sqrt{c^2+|p|^2}-c^2(c^2+|p|^2)^{-1/2}
=2c-(2c^2+|p|^2)(c^2+|p|^2)^{-1/2}.
$$
Notice that $\chi^\prime_2(0)=-|p|<0$ and $\chi^\prime_2(c)$ is continuous.  Suppose $\chi^\prime_2(c^*)=0$ for some $c^*>0$.  Then 
$$
2c^*((c^*)^2+|p|^2)^{1/2}=(2(c^*)^2+|p|^2).
$$
After squaring both sides we observe that this implies $|p|^4=0$, which is a contradiction.  Hence $\chi^\prime_2(c)<0$ for any fixed $|p|>0$ and for all $c\ge 0$.
\end{proof}

The following estimate \eqref{sharpASYMP} will be
important in Section \ref{gnlimit}.

\begin{lemma}
As $c\uparrow\infty$, we have the following pointwise convergence
$$
J\to \mu(v) = (2\pi)^{-N/2} e^{-|v|^2/2}.
$$
For any positive function $h(c) \le A_1 \sqrt{c}$, there is
a uniform constant $A>0$ such that
\begin{equation}
\label{sharpASYMP}
J(h(c)) \le A e^{-h^2(c)/2}.
\end{equation}
By $J(h(c))$, we mean $J(p)$ evaluated at $|p|=h(c)$.
\end{lemma}

The last part of this lemma \eqref{sharpASYMP} is a subtle  asymptotic estimate.  Despite the fact that it is not hard to prove, this estimate is crucial to speed up the convergence in our arguments below for the new rapid convergence rate $1/c^{2-\delta}$ in Theorem \ref{NewtonianLimitThm}.

\begin{proof}  For $0\le x$, we have the following Taylor expansion
\begin{gather}
\sqrt{1+x}=1+\frac{1}{2}x +R(x)
\label{sqrtTAYLOR}
\\
\nonumber
R(x)\eqdef -\frac{1}{8}x^2 (1+x^*)^{-3/2},
\quad
\exists ~ x^*\in [0,x].
\end{gather}
Recall (\ref{juttner}) and note that for $c$ large
$$
K_2(c^2)\sim \left(\frac{\pi}{2}\right)^{1/2}c^{-1} e^{-c^2}+O(c^{-3} e^{-c^2}).
$$
Using this, (\ref{juttner}) and (\ref{sqrtTAYLOR}) we have
\begin{equation*}
J(p)\sim (2\pi)^{-N/2}e^{c^2-cp_0}\to \mu(p).
\end{equation*}
This completes the proof of convergence to a Newtonian Maxwellian.

We now prove the upper bound \eqref{sharpASYMP}.  We use the upper bound in \eqref{jutBOUNDS}.
First the taylor expansion \eqref{sqrtTAYLOR} yields
$$
c^2 - cp_0
=
c^2 - c^2 \sqrt{1+\frac{|p|^2}{c^2}}
=
-
\frac{|p|^2}{2}
-c^2 R\left(\frac{|p|^2}{c^2}\right).
$$
Of course we are working on the set 
$
|p|=h(c) \le A_1 \sqrt{c}.
$
Clearly on this region the expression for $R$ in \eqref{sqrtTAYLOR} gives the bound
$$
\left| c^2 R\left(\frac{|p|^2}{c^2}\right) \right|\le \frac{1}{8}.
$$
Plugging the last two estimates into the upper bound \eqref{jutBOUNDS} yields \eqref{sharpASYMP}. 
\end{proof}

Let $r(c) = \left( \log c \right)^{\alpha_1}$ for $\frac{1}{2} < \alpha_1 <1$.
Given $A>0$, $ \epsilon >0$ there is a uniform  constant $B=B(\epsilon)>0$ such that
\begin{gather}
e^{A  r(c)}
\le 
B c^\epsilon,
\quad
c\ge 1.
\label{slowC}
\end{gather}
Furthermore there is a constant $B_1(k)>0$ such that 
\begin{gather}
e^{-A_1  r^2(c)}
\le 
\frac{B_1(k)}{c^{1+k}},
\quad 
\forall k>0.
\label{fastC}
\end{gather}
These estimates are obvious by after taking the logarithm of both sides.

In the next section we prove 

\section{Global Existence uniformly for $c\ge 1$}\label{geind}

In this section, we prove a uniform (in $c\ge 1$) non-linear estimate which is the key step in our uniform global existence Theorem \ref{ue}.  Global existence of unique mild solutions  can be deduced from this estimate via a standard fixed point argument for small data.  After the proof of Theorem \ref{nonlinLEM} below, we outline how this works.

We split the collision operator \eqref{collisionCM} into its  ``Gain'' and ``Loss'' terms as
\begin{gather*}
\mathcal{Q}_{+}(f,h)
=
\int_{\mathbb{R}^N\times \mathbb{S}^{N-1}} ~ v_c ~ \sigma (g,\theta ) ~ f(p')h(q')~ d\omega dq
\\
\mathcal{Q}_{-}(f,h)
=
\int_{\mathbb{R}^N\times \mathbb{S}^{N-1}} ~ v_c ~ \sigma (g,\theta ) ~ f(p)h(q)~ d\omega dq.
\end{gather*}
With this notation in hand, we prove the following non-linear estimate.

\begin{theorem}
\label{nonlinLEM}
For collision kernels satisfying \eqref{hypS} and \eqref{cutB}, we have 
$$
\int_0^t ~ ds ~ \left| \mathcal{Q}_+^{\#}(f, h)(s,x,p) \right|
\le
A ~\rho_c(x,p)~  \| f^{\#}\|_c \| h^{\#}\|_c .
$$
  Similarly, for the loss term
$$
\int_0^t ~ ds ~ \left| \mathcal{Q}_-^{\#}(f, h)(s,x,p) \right|
\le
A ~\rho_c(x,p)~ \| f^{\#}\|_c \| h^{\#}\|_c .
$$
The constants $A>0$ above are independent of the speed of light.  
\end{theorem}

This estimate may formally look very similar to the Newtonian estimates near Vacuum, e.g. \cite{MR1379589}.  However the main new features of this estimate are twofold.  We have firstly made use of the invariant \eqref{weightCALC} combined with \eqref{cutB} in order to exploit the weak dispersion along characteristics.  Equally important, we have introduced the correct weight \eqref{weightA}  and shown that this allows us to obtain the estimate in a way which is invariant with respect to $c\ge 1$.  Now we are ready for the proof.

\begin{proof} 
We begin with the loss term.  By definition
\begin{multline*}
\int_0^t ds \left| 
\mathcal{Q}_-^{\#}(f,h)(s,x,p) \right|
=
\int_0^t ds
\left| 
\int_{\mathbb{R}^N\times \mathbb{S}^{N-1}} ~ d\omega dq
~ v_c~\sigma_c f(x+s\hat{p},p)h(x+s\hat{p},q)
\right|
\\
=
\int_0^t ds
\left| \int_{\mathbb{R}^N\times \mathbb{S}^{N-1}} ~ d\omega dq~ v_c~\sigma_c ~ f^{\#}(s,x,p) ~ h^{\#}(s,x+s(\hat{p}-\hat{q}),q)
\right|
\\
\le
\| f^{\#}\|_c \| h^{\#}\|_c ~\rho_c(x,p)
\int_0^t ds
 \int_{\mathbb{R}^N}dq~ \int_{\mathbb{S}^{N-1}} d\omega~
 v_c ~\sigma_c (g,\theta) ~ \rho_c(x+s(\hat{p}-\hat{q}),q).
\end{multline*}
We use \eqref{weightA},  Lemma \ref{JSOLbounds} and then \eqref{hypS} to conclude
\begin{multline*}
\int_0^t ds \int_{\mathbb{R}^N}dq~ \int_{\mathbb{S}^{N-1}} d\omega~
v_c ~\sigma(g,\theta) ~ \rho_c(x+s(\hat{p}-\hat{q}),q)
\\
=
\int_0^t ds \int_{\mathbb{R}^N}dq~ \int_{\mathbb{S}^{N-1}} d\omega~v_c ~\sigma (g,\theta )~ J^{\beta}(q)
e^{- \alpha q_0 |x+s(\hat{p}-\hat{q})|^2/c}
\\
\le
\int_0^t ds \int_{\mathbb{R}^N}dq~ \int_{\mathbb{S}^{N-1}} d\omega~v_c ~\sigma (g,\theta ) ~e^{-\beta|q|}
e^{- \alpha  |x+s(\hat{p}-\hat{q})|^2}
\\
=
 \int_{\mathbb{R}^N}dq ~ \int_{\mathbb{S}^{N-1}} d\omega ~ v_c ~ \sigma (g,\theta ) ~ e^{-\beta|q|} ~
\int_0^t ds ~e^{- \alpha  |x+s(\hat{p}-\hat{q})|^2}
\\
\le
 \int_{\mathbb{R}^N}dq ~ \int_{\mathbb{S}^{N-1}} d\omega ~ \sigma (g,\theta ) ~ e^{-\beta|q|}
 \le
 A.
\end{multline*}
In getting to the last line above we utilized the standard estimate
$$
\int_0^t ds ~e^{- \alpha  |x+s(\hat{p}-\hat{q})|^2}
\le
\sqrt{\frac{\pi}{\alpha}}\frac{1}{|\hat{p}-\hat{q}|}.
$$
This follows easily by completing the square in the exponent of the exponential; 
see for instance \cite{MR1379589}. We have also used in the same spot the estimate
$$
\frac{v_c}{|\hat{p}-\hat{q}|} \le
\frac{c}{2}\left| \frac{p}{p_0}-\frac{q}{q_0}\right| \frac{1}{|\hat{p}-\hat{q}|}
=
 \frac{1}{2},
$$
which follows from \eqref{moller}.
In our estimates for the loss term we do not need to use the restriction into the set 
\eqref{cutB} at all.  

Furthermore, we explain the last inequality which will follow from \eqref{hypS}
\begin{gather*}
 \int_{\mathbb{R}^N}dq ~ \int_{\mathbb{S}^{N-1}} d\omega ~ \sigma (g,\theta ) ~ e^{-\beta|q|}
 \le
 A\int_{\mathbb{R}^N}dq ~ \int_{\mathbb{S}^{N-1}} d\omega ~ \left( 1+ ~ g^{-\gamma}\right)  ~ e^{-\beta|q|}
 \le
 A.
\end{gather*}
The last bound follows easily from our independent of $c$ bound on $g$ in Lemma \ref{lowerBB}.

Next we consider the gain term, for which  we have the estimate
\begin{multline*}
\int_0^t ds \left| 
\mathcal{Q}_+^{\#}(f,h)(s,x,p) \right|
\\
=
\int_0^t ds
\left| \int_{\mathbb{R}^N \times \mathbb{S}^{N-1}}dq~ d\omega~ v_c~\sigma~f(s,x+s\hat{p},p')h(s,x+s\hat{p},q')
\right|
\\
=
\int_0^t ds
\left| \int_{\mathbb{R}^N}dq~ \int_{\mathbb{S}^{N-1}} d\omega~ v_c\sigma  ~ f^{\#}(s,x+s(\hat{p}-\hat{p}'),p') ~ h^{\#}(s,x+s(\hat{p}-\hat{q}'),q')
\right|
\\
\le
 \| f^{\#}\|_c \| h^{\#}\|_c 
\int_0^t ds
 \int_{\mathbb{R}^N \times \mathbb{S}^{N-1}}dq~ d\omega~ v_c ~ \sigma ~ \rho_c(x+s(\hat{p}-\hat{p}'),p')\rho_c(x+s(\hat{p}-\hat{q}'),q').
\end{multline*}
Next we use the weights $\rho_c$ and their arguments to control the remaining integrations.

Using \eqref{weightA}, \eqref{weightCALC} and \eqref{collisionalCONSERVATION}
we have 
\begin{gather*}
\rho_c(x+s(\hat{p}-\hat{p}'),p')
\rho_c(x+s(\hat{p}-\hat{q}'),q')
=
\rho_c(x,p)\rho_c(x+s(\hat{p}-\hat{q}),q) e^{-\delta t^2}.
\end{gather*}
On $\mathcal{B}_c$, in \eqref{cutB}, we have
$$
\delta = c^3 \left(\frac{1}{p_0}+\frac{1}{q_0}  - \frac{1}{p_0'} - \frac{1}{q_0'} \right) \ge - h_c.
$$
Recalling \eqref{cutHC}, we have the bound
$$
e^{\delta t^2} 
\ge 
e^{ -  h_c t^2}
= e^{ - B - a \alpha q_0 \left|x+t\left( \hat{p} - \hat{q}\right)\right|^2/c}.
$$
Combining this estimate with the weight \eqref{weightA} and recalling that $0<a<1$ we obtain the upper bound
\begin{gather*}
\rho_c(x+s(\hat{p}-\hat{q}),q) e^{-\delta t^2}
\le
 e^{  B} e^{-(1- a) \alpha q_0 \left|x+t\left( \hat{p} - \hat{q}\right)\right|^2/c} ~ J^{\beta}(q).
\end{gather*}
We conclude  
\begin{multline*}
\int_0^t ds \left| 
\mathcal{Q}_+^{\#}(f,g)(s,x,p) \right|
\\
\le
\| f^{\#}\|_c \| h^{\#}\|_c  ~ \rho_c(x,p) ~ 
\int_0^t ds
 \int_{\mathbb{R}^N \times \mathbb{S}^{N-1}} dq d\omega~v_c~ \sigma_c ~  e^{B-(1- a) \alpha q_0 \left|x+t\left( \hat{p} - \hat{q}\right)\right|^2/c} ~ J^{\beta}(q)
 \\
 \le
A~  \| f^{\#}\| \| h^{\#}\| ~ \rho_c(x,p).
\end{multline*}
The final integral is controlled in the same way as for the loss term.
\end{proof}

Now that our  non-linear estimate has been established, it becomes a standard exercise in well known and well exposited techniques, e.g. \cite{MR760333,MR0475532} and for instance \cite{MR1379589}, 
to prove global existence of unique positive mild solutions.
In particular  
the fixed point argument and Kaniel-Shinbrot iteration scheme described in
 \cite{MR2217287}
can be 
used without any modification 
in combination with our scale invariant non-linear estimate in Theorem \ref{nonlinLEM}
in order to prove global in time existence of unique mild solutions 
with a small data condition which is uniform in the speed of light $c\ge 1$.

In the next section we consider the Newtonian Limit.

\section{Global Newtonian Limit as $c\to\infty$}\label{gnlimit}

In this section, we work with the relativistic Boltzmann equation \eqref{RBF} in the forward mild form
\begin{gather}
f_c(t,x,p) = S_{c}(t)f_{0,c}(x,p) + \int_0^t ~ ds ~ S_{c}(t-s)\mathcal{Q}_c(s, x, p).
\label{mildSOL}
\end{gather}
We are using the definition $S_{c}(t)f(x,p) \eqdef f(x- t \hat{p},p)$, where we recall from \eqref{normV} that
$
\hat{p}  =c\frac{p}{p_0}=\frac{p}{\sqrt{1+|p|^2/c}}
$
clearly depends upon $c$. 
Notice that this mild formulation can be easily interchanged with the backward formulation that we used earlier \eqref{mildRBE}; this one is however more convenient for the Newtonian limit.
Existence of such solutions is covered by Theorem \ref{ue}.  
Above we also use the abbrieviated notation $\mathcal{Q}_c = \mathcal{Q}_c(f_c, f_c)$.

We also study the Newtonian Boltzmann equation \eqref{RBFnewt}
 in the forward mild form
\begin{gather}
f(t,x,p) = S(t)f_{0}(x,p) + \int_0^t ~ ds ~ S(t-s)Q(s, x, p).
\label{mildSOLnewt}
\end{gather}
We are here also using the definition $S(t)f(x,p) \eqdef f(x- t p,p)$.  The existence theory for \eqref{mildSOLnewt} was explained in Section \ref{mainres:sec}.  For the limit system, we have continuity.

\begin{lemma}\label{controlDIFF}  Let  the conditions in Theorem \ref{NewtonianLimitThm} hold.   
Given $T >0$ and any  $0<\delta < 1$ 
there is a constant  $A = A(\delta, T)>0$ such that for any $0\le t \le T$ we have
\begin{equation}
\sup_{|u| \le |h|} 
\left(\| \tau_  u^x f(t) - f(t)  \|_{L^1_p L^\infty_x}
+
\| \tau_  u^p f(t)  - f(t)  \|_{L^1_p L^\infty_x}
\right)
\le
A  |h|^{\delta},
\quad 
|h| < 1.
\notag
\end{equation}
\end{lemma}

Lemma \ref{controlDIFF} proves the continuity of the mild form of the Newtonian Boltzmann equation which will be an important part of our convergence theorem.

\begin{proof}  To estimate the $x$ difference  we subtract the two  equations \eqref{mildSOLnewt} as
\begin{eqnarray}
\tau_  h^x f(p) - f(p)
&=&
\tau_  h^x S(t)f_{0}(x, p) - S(t) f_{0}(x, p)
\nonumber
\\
&&
+ \int_0^t ~ ds ~ \tau_  h^x S(t-s) Q(s, x, p)
 - 
 S(t-s) Q(s, x, p), 
 \nonumber
\end{eqnarray}
which can then be further expanded as
\begin{gather}
 =
\tau_  h^x S(t)f_{0}(x, p) - S(t) f_{0}(x, p)
 \label{Lt1}
\\
+ \int_0^t  ~  \int d\omega dq ~ \mathcal{K}_\infty ~
\left( \tau_  h^x S(t-s) f(\bar{p}^{\prime}) - S(t-s) f(\bar{p}^{\prime})\right) ~ \tau_  h^x S(t-s) f(\bar{q}^{\prime})
  \label{Lt3}
  \\
+ \int_0^t ~ ds ~  \int_{\mathbb{R}^N\times \mathbb{S}^{N-1}} d\omega dq ~ \mathcal{K}_\infty ~
S(t-s) f(\bar{p}^{\prime}) ~ \left(\tau_  h^x S(t-s) f(\bar{q}^{\prime}) - S(t-s) f(\bar{q}^{\prime}) \right)
  \nonumber
 \\
- \int_0^t ~ ds ~ \int_{\mathbb{R}^N\times \mathbb{S}^{N-1}} d\omega dq ~ \mathcal{K}_\infty  ~
\left( \tau_  h^x S(t-s) f(p) - S(t-s) f(p)\right) ~ \tau_  h^x S(t-s) f(q)
  \nonumber
  \\
- \int_0^t ~ ds ~  \int_{\mathbb{R}^N\times \mathbb{S}^{N-1}} d\omega dq ~ \mathcal{K}_\infty ~
S(t-s) f(p) ~ \left(\tau_  h^x S(t-s) f(q) - S(t-s) f(q) \right).
  \nonumber
\end{gather}
The $L^1_p L^\infty_x$ norm of  \eqref{Lt1} is easily controlled 
\begin{gather*}
\| \tau_  h^x S(t)f_{0} - S(t) f_{0} \|_{L^1_p L^\infty_x}
=
\| \tau_  h^x f_{0} -  f_{0} \|_{L^1_p L^\infty_x},
\end{gather*}
which is simply a statement of translation invariance of the $L^\infty_x$ norm.

For the other terms, we split into an unbounded and a bounded region
\begin{equation}
1 = {\bf 1}_{|p|+|q| \ge R}
+
{\bf 1}_{|p|+|q| < R},
 \label{largeSMALLsplit}
\end{equation}
where we fix a suitable $R>0$ below.  

We take the supremum in $x$ and integrate over $p$.  On the unbounded region, $|p|+|q|>R$, we see that each of the terms  \eqref{Lt3} and after  are bounded above by
\begin{equation}
b_2 \int_0^t ~ ds ~ \int_{\mathbb{R}^N\times\mathbb{R}^N\times \mathbb{S}^{N-1}} d\omega dq dp 
~ \mathcal{K}_\infty  ~  
{\bf 1}_{|p|+|q| \ge R} ~ \mu^\beta(p) ~ \mu^\beta(q)
 \le b_2 \sigma_0 T  e^{-\beta R^2/4}.
 \label{largeBDD}
\end{equation}
We have used the 
uniform bound for solutions in Theorem \ref{NewtonianLimitThm}, \eqref{weightINF}, and \eqref{kernelUPPERbds}.

On the bounded region we have the following upper bound for the term \eqref{Lt3}
\begin{multline*}
\le 
  \int_0^t  ds ~ \int_{\mathbb{R}^N\times\mathbb{R}^N\times \mathbb{S}^{N-1}} d\omega dq dp ~ {\bf 1}_{|p|+|q| < R} ~ 
  \mathcal{K}_\infty  ~
\| \tau_  h^x f(\bar{p}^{\prime}) -   f(\bar{p}^{\prime}) \|_{L^\infty_x}  \|  \tau_  h^x f(\bar{q}^{\prime})  \|_{L^\infty_x}
\\
\le 
 \int_0^t  ds ~ \int_{\mathbb{R}^N\times\mathbb{R}^N\times \mathbb{S}^{N-1}} 
 d\omega dq dp   ~
  {\bf 1}_{|p|+|q| < 4R}~
  \mathcal{K}_\infty
\| \tau_  h^x f(p) -   f(p) \|_{L^\infty_x}  \|  \tau_  h^x f(q)  \|_{L^\infty_x}
\\
\le
AR \sigma_0
 \int_0^t  ds ~ \int_{\mathbb{R}^N\times\mathbb{R}^N\times \mathbb{S}^{N-1}} d\omega dq dp ~
\| \tau_  h^x f(p) -   f(p) \|_{L^\infty_x}  \|\tau_  h^x  f(q)  \|_{L^\infty_x}
\\
\le A R   \int_0^t  ds ~ \| \tau_  h^x f -   f \|_{L^1_p L^\infty_x}  \|  f  \|_{L^1_p L^\infty_x}
\\
\le A R   \int_0^t  ds ~ \| \tau_  h^x f -   f \|_{L^1_p L^\infty_x}  \|  f^{\#}  \|_\infty
\le  A b_2 R   \int_0^t  ds ~ \| \tau_  h^x f -   f \|_{L^1_p L^\infty_x}.  
\end{multline*}
Above the constant $A>0$ is independent of  $R$. 
We  just used the pre-post collisional change of variables, \eqref{kernelUPPERbds} and the assumption of Theorem \ref{NewtonianLimitThm} for the uniform bound in terms of $b_2$.  
Notice that all the other terms  after \eqref{Lt3} can be treated in exactly the same way.  We have shown
\begin{multline*}
\| \tau_  h^x f - f \|_{L^1_p L^\infty_x}
\le
 \| \tau_  h^x f_{0} - f_{0} \|_{L^1_p L^\infty_x}
 +
A_1 T e^{-\beta R^2/4}
\\
+
A_2  R   \int_0^t  ds ~ \| \tau_  h^x f -   f \|_{L^1_p L^\infty_x}.  
\end{multline*}
By the Gronwall inequality 
\begin{gather*}
\| \tau_  h^x f - f \|_{L^1_p L^\infty_x}
\le
\left(
 \| \tau_  h^x f_{0} - f_{0} \|_{L^1_p L^\infty_x}
 +
A_1 T e^{-\beta R^2/4}
\right)
e^{A_2  R T}.
\end{gather*}
Now choose $R = \left( - \log |h| \right)^{\alpha_1}$ for any $\frac{1}{2} < \alpha_1 <1$.  Using 
\eqref{slowC} and \eqref{fastC} with $c = 1/|h|$ and $|h| < 1$
we obtain for any small $\epsilon >0$ that
\begin{gather*}
\| \tau_  h^x f - f \|_{L^1_p L^\infty_x}
\le
A \left(
|h|^{-\epsilon} \| \tau_  h^x f_{0} - f_{0} \|_{L^1_p L^\infty_x}
 +
|h|
\right).
\end{gather*}
The proof of the estimate for $\| \tau_  h^x f - f \|_{L^1_p L^\infty_x}$ is completed by using Lipschitz property on the initial data in Theorem \ref{NewtonianLimitThm}.

For the estimate of $\| \tau_  h^p f - f \|_{L^1_p L^\infty_x}$, we write out the difference of \eqref{mildSOLnewt} as
\begin{multline}
\tau_  h^p f(p) - f(p)
=
\tau_  h^p S(t)f_{0}(x, p) - S(t) f_{0}(x, p)
\\
+ \int_0^t ~ ds ~ \tau_  h^p S(t-s) Q(s, x, p)
 - 
 S(t-s) Q(s, x, p), 
 \\
 = T_1 + T_2.
\notag
\end{multline}
We can expand the first term as
\begin{multline}
T_1
=
f_{0}(x-t(p + h), p+h) -  f_{0}(x-tp, p)
\\
=
\left\{ f_{0}(x-t(p + h), p+h) -  f_{0}(x-t(p + h), p)  \right\}
\\
+
\left\{ f_{0}(x-t(p + h), p) -  f_{0}(x-tp, p)  \right\}.
\notag
\end{multline}
 Thus,
\begin{multline}
\|  T_1  \|_{L^1_p L^\infty_x}
\le 
\|  \tau_  h^p f_{0} - f_{0} \|_{L^1_p L^\infty_x}
+
\|  \tau_ {th}^x f_{0} - f_{0} \|_{L^1_p L^\infty_x}
\\
\le 
\|  \tau_  h^p f_{0} - f_{0} \|_{L^1_p L^\infty_x}
+
T\|  \tau_ {h}^x f_{0} - f_{0} \|_{L^1_p L^\infty_x}.
\notag
\end{multline}
Here we used the translation invariance of the norm to pull the fixed constant out of the supremum.
More precisely, for any integer $n\ge 1$ we used the estimate
\begin{multline}
\|  \tau_ {nh}^x f_{0} - f_{0} \|_{L^1_p L^\infty_x}
\le 
\sum_{i= 0}^n
\|  \tau_ {ih}^x f_{0} - \tau_ {(i-1)h}^xf_{0} \|_{L^1_p L^\infty_x}
\\
\le 
\sum_{i= 0}^n
\|  \tau_ {h}^x f_{0} - f_{0} \|_{L^1_p L^\infty_x}
\le 
n
\|  \tau_ {h}^x f_{0} - f_{0} \|_{L^1_p L^\infty_x},
\notag
\end{multline}
which follows because of the translation invariance of these norms.  Technically, we should have supposed that $T$ was an integer, but this does not cause difficulties.

The other term $T_2$ involves the collision operator.   In this spot we  use the Newtonian variables \eqref{omegaREP}.
We further split as $T_2 = T_2' + T_3$ with
 $T_2'$ given by
\begin{multline}
T_2'  =
 \int_0^t  ds ~  \int d\omega dq ~ \left(\tau_h^p\mathcal{K}_\infty\right)(p,q,\omega) ~  f(s, \hat{x}_h, \bar{p}_h') ~ 
 f(s, \hat{x}_h, \bar{q}_h')
\\
-  \int_0^t  ds ~  \int d\omega dq ~ \mathcal{K}_\infty(p,q,\omega) ~  f(s, \hat{x}, \bar{p}') ~ f(s,\hat{x}, \bar{q}'). 
\notag
\end{multline}
The remaining term $T_3$ is then given by
\begin{multline*}
T_3  =
 -\int_0^t  ds ~  \int d\omega dq    ~
 \left(\tau_h^p\mathcal{K}_\infty\right)(p,q,\omega) ~
  f(s, \hat{x}_h, p+h) ~ f(s, \hat{x}_h, q)
  \\
 +\int_0^t  ds ~  \int d\omega dq    ~
\mathcal{K}_\infty (p,q,\omega)~ 
f(s,\hat{x}, p) ~ f(s, \hat{x}, q).
\end{multline*}
We have used the following simplifying notation:
 $\hat{x}_h = x-(t-s)(p + h)$ and similarly $\hat{x} = x-(t-s) p$.
Also with   \eqref{omegaREP}
we define 
$$
\bar{p}'_h = p+h+\omega \cdot \left(q-p-h\right)\omega,
\quad 
\bar{q}^{\prime }_h =q-\omega \cdot \left(q-p-h\right)\omega.
$$ 
For both $T_2'$ and $T_3$, we use the splitting \eqref{largeSMALLsplit}.  On the large region ${\bf 1}_{|p|+|q| \ge R}$, we again have \eqref{largeBDD} for both $T_2'$ and $T_3$ (the small $|h| < 1$ causes no disruption).

On the bounded region,
 ${\bf 1}_{|p|+|q| \le R}$,
we expand $T_2'$ as
\begin{multline}
T_2'=
 \int_0^t  ds ~  \int d\omega dq ~{\bf 1}_{|p|+|q| \le R}  \left\{ \left(\tau_h^p\mathcal{K}_\infty\right) - \mathcal{K}_\infty  \right\}~   f(s,\hat{x}_h, \bar{p}_h') ~ f(s, \hat{x}_h, 
\bar{q}'_h)
 \\
 +
 \int_0^t  ds ~  \int d\omega dq ~{\bf 1}_{|p|+|q| \le R}~  \mathcal{K_\infty} ~ \left\{ f(s, \hat{x}_h, \bar{p}_h') 
- f(s, \hat{x}_h, \bar{p}')  \right\} f(s, \hat{x}_h, \bar{q}_h')
\\
+
 \int_0^t  ds ~  \int d\omega dq ~{\bf 1}_{|p|+|q| \le R}~ \mathcal{K_\infty} ~  f(s, \hat{x}_h, \bar{p}')   \left\{f(s, \hat{x}_h, \bar{q}_h') - f(s, \hat{x}_h, \bar{q}') \right\}
\\
+
 \int_0^t  ds ~  \int d\omega dq ~ {\bf 1}_{|p|+|q| \le R}~ \mathcal{K_\infty} ~  f(s, \hat{x}_h, \bar{p}')   \left\{f(s, \hat{x}_h, \bar{q}') - f(s, \hat{x}, \bar{q}') \right\}
\\
+
 \int_0^t  ds ~  \int d\omega dq ~ {\bf 1}_{|p|+|q| \le R}~ \mathcal{K_\infty} ~ f(s, \hat{x}, \bar{q}')   \left\{f(s, \hat{x}_h, \bar{p}') - f(s, \hat{x}, \bar{p}') \right\}.
\notag
\end{multline}
As a result of the hypothesis in Theorem \ref{NewtonianLimitThm},  since $\mathcal{K}_\infty$ is Lipschitz continuous  we have
$$
\left| \left(\tau_h^p\mathcal{K}_\infty\right) - \mathcal{K}_\infty  \right|
e^{-\beta (|p|^2+|q|^2)/4}
\le
A  |h|
e^{-\beta (|p|^2+|q|^2)/8}.
$$
Here there could be a smooth polynomial appearing in $|p|$ and/or $|q|$ which is controlled by the exponential decay.
With this estimate, the first term in $T_2'$ integrated over $dp$, has the following upper bound estimate
\begin{multline}
 \int_0^t  ds ~  \int d\omega dq dp ~ \left| \left(\tau_h^p\mathcal{K}_\infty\right) - \mathcal{K}_\infty  \right|~  
 \| f(\bar{p}'_h) \|_{L^\infty_x} ~  \| f(\bar{q}'_h) \|_{L^\infty_x} 
 \\
 \le
A ~ b_1^2 ~ T \int d\omega dq dp ~ \left| \left(\tau_h^p\mathcal{K}\right) - \mathcal{K}  \right|~  
e^{-\beta (|p|^2+|q|^2)/4}
  \le
A  b_1^2  T  |h|.
\notag
\end{multline}
The estimate above follows from the 
uniform bound for solutions in Theorem \ref{NewtonianLimitThm} and \eqref{weightINF} as well as the pre-post invariance of the exponentials.  Notice that the translation of a small $h$ does not cause difficulties.

For the next terms in $T_2'$, with the bounded region, we will use several times that
$$
|\bar{q}'_h - \bar{q}'|+|\bar{p}'_h - \bar{p}'| \le 4 |h|.
$$
This is a trivial consquence of \eqref{omegaREP}.
Now the estimates for the next terms in $T_2'$ work as in a previous case.  In particular
\begin{multline}
 \int_0^t  ds ~  \int d\omega dq dp ~ {\bf 1}_{|p|+|q| < R} ~
  \left|\mathcal{K}_\infty\right| ~ 
  \| f(\bar{p}_h') - f(\bar{p}')  \|_{L^\infty_x}  \| f(\bar{q}_h')\|_{L^\infty_x}
\\
\le A'  \int_0^t  ds ~    \int d\omega dq dp ~ {\bf 1}_{|p|+|q| < R} ~ 
 \left|\mathcal{K}_\infty\right| ~ 
  \sup_{|u| \le 4 |h|}\| (\tau^p_u f)(\bar{p}') - f(\bar{p}')  \|_{L^\infty_x}  
  \mu^{\beta/2} (\bar{q}')
  \\
\le A'  \int_0^t  ds ~    \int d\omega dq dp ~ {\bf 1}_{|p|+|q| < 4R} ~ 
 \left|\mathcal{K}_\infty\right| ~ 
  \sup_{|u| \le 4 |h|}\| (\tau^p_u f)(p) - f(p)  \|_{L^\infty_x}  
  \mu^{\beta/2} (q)
    \\
\le A'  R  \int_0^t  ds ~   
  \sup_{|u| \le 4 |h|}\| \tau^p_u f - f  \|_{L^1_pL^\infty_x}  
  \\
  \le 4 A'   R  \int_0^t  ds ~   
  \sup_{|u| \le  |h|}\| \tau^p_u f - f  \|_{L^1_pL^\infty_x}.  
\notag
\end{multline}
Here we have used the estimate for solutions in Theorem \ref{NewtonianLimitThm} and \eqref{kernelUPPERbds}.
All of the other remaining estimates in  $T_2'$ are exactly the same as this one, or a previous estimate in the $x$ variables given earlier in this proof.  The loss term estimates are the same but easier because they do not require a pre-post change of variable.

For the $T_3$ term we employ the following splitting
\begin{multline*}
-T_3  
=
 \int_0^t  ds ~  \int d\omega dq ~ 
 \left\{ \left(\tau_h^p\mathcal{K}_\infty\right) - \mathcal{K}_\infty  \right\} ~  f(s, \hat{x}_h, p+h) f(s,\hat{x}_h, q)
 \\
 +
 \int_0^t  ds ~  \int d\omega dq ~ \mathcal{K}_\infty ~  \left\{ f(s, \hat{x}_h, p+h) ~ 
-  f(s,\hat{x}_h, p)  \right\} f(s, \hat{x}_h, q)
\\
+
 \int_0^t  ds ~  \int d\omega dq ~ \mathcal{K}_\infty ~  \left\{ f(s, \hat{x}_h, q) ~ 
-  f(s, \hat{x}, q) \right\} f(s,\hat{x}_h, p) 
\\
+
 \int_0^t  ds ~  \int d\omega dq ~  \mathcal{K}_\infty ~  \left\{ f(s, \hat{x}_h, p) ~ 
-  f(s, \hat{x}, p) \right\} f(s,\hat{x}, q).
\end{multline*}
All of these terms have a corresponding expression in the $T_2'$ expansion above (except there are no primed variables here, which simplifies things).  Thus the estimates for $T_3 $ are exactly the same as the ones for $T_2'$, and they are proven in the same way.

Collecting  all of these estimates, we have shown
\begin{multline*}
\| \tau_  h^p f - f \|_{L^1_p L^\infty_x}
\le
\|  \tau_  h^p f_{0} - f_{0} \|_{L^1_p L^\infty_x}
+
A' T\sup_{|u| \le  |h| }\|  \tau_ u^x f_{0} - f_{0} \|_{L^1_p L^\infty_x}
 +
A_1 T e^{-\beta R^2/2}
\\
+
A T |h|
+
A_2'  R   \int_0^t  ds ~
\sup_{|u| \le  |h| }\left( \| \tau_  h^p f -   f \|_{L^1_p L^\infty_x}
+ \|  \tau_ u^x f - f \|_{L^1_p L^\infty_x} \right).  
\end{multline*}
We apply the Lipschitz condition on the initial data in Theorem \ref{NewtonianLimitThm}
as well as the prior estimate for $\| \tau_ h^x f - f \|_{L^1_p L^\infty_x}$   to obtain
for any small $\epsilon >0$ that
\begin{multline*}
\sup_{|u| \le  |h| }
\| \tau_  u^p f - f \|_{L^1_p L^\infty_x}
\le
A(T) \left(  |h|^{1-\epsilon}  + e^{ - \beta R^2 /4} \right)
\\
+
A'  R   \int_0^t  ds ~
\sup_{|u| \le  |h| } \| \tau_u^p f -   f \|_{L^1_p L^\infty_x}.
\end{multline*}
By the Gronwall inequality 
\begin{gather*}
\sup_{|u| \le  |h| }
 \| \tau_  u^p f - f \|_{L^1_p L^\infty_x}
\le
A(T) \left(  |h|^{1-\epsilon}  + e^{ - \beta R^2 /4} \right)
e^{A'  R T}.
\end{gather*}
As before choose $R = \left( - \log |h| \right)^{\alpha_1}$ for any $\frac{1}{2} < \alpha_1 <1$. 
 Using 
\eqref{slowC} and \eqref{fastC} with $c = 1/|h|$ and $|h| < 1$
we obtain for any small $\delta >0$ the statement of Lemma \ref{controlDIFF}.
\end{proof}

We now give our main estimate for the difference of the two collision operators.

\begin{lemma}  
\label{nonLINasy}
Assume the conditions of Theorem \ref{NewtonianLimitThm} hold.  
Suppose that $ c \ge 1$, $\frac{1}{2}<\alpha_1<1$, and $0\le t \le T$.  Then for any  $\varepsilon\in(0,2]$ we have
\begin{gather*}
\| S_{c}(t)\mathcal{Q}_{c}(t) -  S(t)\mathcal{Q}(t)\|_{L^1_p L^\infty_x}
\le
A \left( \log c \right)^{\alpha_1}
\| f_c - f \|_{L^1_p L^\infty_x} 
+
\frac{ A(T)}{c^{2-\varepsilon}}.
\end{gather*}
These collision operators are defined as in \eqref{collisionGS} and \eqref{RBFnewt}.
Furthermore, as usual, the positive constants $A$, $A(T)$ are uniform in $c$.
\end{lemma}

\begin{proof}  From \eqref{collisionGS} and \eqref{RBFnewt} we have 
\begin{multline}
S_{c}(t)\mathcal{Q}_{c}
-
S(t)\mathcal{Q}
= 
\int_{\mathbb{R}^N\times \mathbb{S}^{N-1}} ~d\omega dq ~ \mathcal{K}_{c} ~ S_{c}(t)f_{c}(p^{\prime }) S_{c}(t)f_{c}(q^{\prime})
\\
-\int_{\mathbb{R}^N\times \mathbb{S}^{N-1}} ~d\omega dq ~ \mathcal{K}_{c}  ~ S_{c}(t)f_{c}(p) S_{c}(t)f_{c}(q)
\\
-\int_{\mathbb{R}^N\times \mathbb{S}^{N-1}} ~d\omega dq ~ \mathcal{K}_{\infty}  ~ S(t)f(\bar{p}^{\prime }) S(t)f(\bar{q}^{\prime})
\\
+\int_{\mathbb{R}^N\times \mathbb{S}^{N-1}} ~d\omega dq ~ \mathcal{K}_{\infty}  ~ S(t) f(p) S(t) f(q).
\label{termsz}
\end{multline}
As usual, we split the integrals via 
$$
1={\bf 1}_{|p|+|q| \le r(c)} + {\bf 1}_{|p|+|q| > r(c)}.
$$
Here $r(c) = \left(\log c\right)^{\alpha_1}$ with $1/2 < \alpha_1 <1$.  

When the momentum are large we have the estimate
\begin{multline}
\label{extraCest}
\left| 
\int_{\mathbb{R}^N\times \mathbb{S}^{N-1}} ~d\omega dq dp ~ {\bf 1}_{|p|+|q| > r(c)} ~ \mathcal{K}_{c}  ~ S_{c}(t)f_{c}(p^{\prime }) S_{c}(t)f_{c}(q^{\prime})
\right|
\\
\le b_2^2
\int_{\mathbb{R}^N\times \mathbb{S}^{N-1}} ~d\omega dq dp ~ {\bf 1}_{|p|+|q| > r(c)} ~ \mathcal{K}_{c}  ~ J^{\beta}(p') J^{\beta}(q')
\\
= b_2^2
\int_{\mathbb{R}^N\times \mathbb{S}^{N-1}} ~d\omega dq dp ~ {\bf 1}_{|p|+|q| > r(c)} ~ \mathcal{K}_{c}  ~ J^{\beta}(p) J^{\beta}(q)
\\
\le 
b_2^2
J^{\beta/2}(r(c)) 
\int_{\mathbb{R}^N\times \mathbb{S}^{N-1}} ~d\omega dq dp ~ {\bf 1}_{|p|+|q| > r(c)} ~ \mathcal{K}_{c}  ~ J^{\beta/2}(p) J^{\beta/2}(q)
\\
\le 
A b_2^2 \sigma_0 J^{\beta/2}(r(c)) 
\int_{\mathbb{R}^N\times \mathbb{S}^{N-1}} ~d\omega dq dp ~ {\bf 1}_{|p|+|q| > r(c)}   ~ J^{\beta/4}(p) J^{\beta/4}(q)
\\
\le 
A  J^{\beta/2}(r(c)) 
\\
\le
 A e^{-\beta (r(c))^2/4}
 \le 
 \frac{A}{c^{2}}. 
\end{multline}
Above we have used the estimate for $\mathcal{K}_{c}(p,q,\omega)$ from Lemma \ref{asymptoticBB}.
It is crucial that we gain back the square in the relativistic Maxwellian as in 
\eqref{sharpASYMP}, in order to use the estimate \eqref{fastC} in the last line.
We have also used the upper bound of $J^{\beta}(p)$ in the weight of the norm $\| \cdot \|_c$ from Theorem \ref{NewtonianLimitThm}, and the pre-post collisional invariance \eqref{collisionalCONSERVATION}.  The independent of the speed of light bound for $J(p)$ in Lemma \ref{JSOLbounds} would not be enough.  All the other terms on the r.h.s. of \eqref{termsz} can be easily controlled by the same upper bound when the momentum are large.

For bounded momentum we use more splitting.  Notice that the momentum post-collision depend upon the speed of light.  We expand the terms on the r.h.s. of \eqref{termsz}:
\begin{multline}
\int ~d\omega  dq  ~  {\bf 1}_{|p|+|q| \le r(c)} ~ \mathcal{K}_{c}  ~ 
S_{c}(t)f_{c}(p^{\prime}) S_{c}(t)f_{c}(q^{\prime})
\nonumber
\\
-
\int  ~d\omega  dq    ~  {\bf 1}_{|p|+|q| \le r(c)} ~ \mathcal{K}_{c}  ~ S_{c}(t)f_{c}(p) S_{c}(t)f_{c}(q)
\nonumber
\\
-\int  ~ d\omega  dq     ~  {\bf 1}_{|p|+|q| \le r(c)}~ \mathcal{K}_{\infty}  ~ 
S(t)f(\bar{p}^{\prime })S(t)f(\bar{q}^{\prime})
\nonumber
\\
+\int  ~ d\omega  dq     ~  {\bf 1}_{|p|+|q| \le r(c)} ~ \mathcal{K}_{\infty}  ~ 
S(t) f(p)S(t) f(q)
\\
= 
\mathbb{I}
+
\mathbb{II}
+
\mathbb{III}.
\nonumber
\end{multline}
For the first term we have  
\begin{gather*}
\mathbb{I}
= 
\int  ~ d\omega dq     ~ 
 {\bf 1}_{|p|+|q| \le r(c)}
~ \left\{ \mathcal{K}_{c}  -   \mathcal{K}_{\infty} \right\}  ~ 
\left\{
S(t) f(\bar{p}^{\prime}) S(t) f(\bar{q}^{\prime})
-
S(t) f(p) S(t) f(q)
\right\}.
\end{gather*}
For the second term in our splitting above is
\begin{gather}
\mathbb{II}
= 
\int  ~ d\omega  dq     ~ 
 {\bf 1}_{|p|+|q| \le r(c)}
~  \mathcal{K}_{c} ~ [S_{c}(t) f_{c}(p^{\prime})- S_{c}(t) f(p^{\prime})]S_{c}(t)f_{c}(q^{\prime})
\label{t9}
\\
+
\int  ~ d\omega  dq     ~ 
 {\bf 1}_{|p|+|q| \le r(c)}
~ \mathcal{K}_{c} ~ [S_{c}(t) f(p^{\prime})- S(t) f(p^{\prime})]S_{c}(t)f_{c}(q^{\prime})
\label{t9a}
\\
+
\int  ~ d\omega  dq     ~ 
 {\bf 1}_{|p|+|q| \le r(c)}
~ \mathcal{K}_{c}~ [S_{c}(t) f_{c}(q^{\prime})-S_{c}(t) f(q^{\prime})] S(t) f(p^{\prime})
\label{t10a}
\\
+
\int ~ d\omega  dq     ~ 
 {\bf 1}_{|p|+|q| \le r(c)}
~ \mathcal{K}_{c} ~ [S_{c}(t) f(q^{\prime})-S(t) f(q^{\prime})] S(t) f(p^{\prime})
\label{t11}
\\
+
\int ~ d\omega  dq     ~ 
 {\bf 1}_{|p|+|q| \le r(c)}
~ \mathcal{K}_{c} ~ [S(t) f(q^{\prime})-S(t) f(\bar{q}^{\prime})] S(t) f(p^{\prime})
\label{t12}
\\
+
\int  ~ d\omega  dq     ~ 
 {\bf 1}_{|p|+|q| \le r(c)}
~ \mathcal{K}_{c} ~ [S(t) f(p^{\prime})- S(t) f(\bar{p}^{\prime})]S(t)f(\bar{q}^{\prime}).
\label{t10}
\end{gather}
The third term above can be expanded as
\begin{multline*}
-\mathbb{III}
= 
\int  ~ d\omega  dq     ~ 
 {\bf 1}_{|p|+|q| \le r(c)}
~  \mathcal{K}_{c} ~ [S_{c}(t) f_{c}(p)- S_{c}(t) f(p)]S_{c}(t)f_{c}(q)
\\
+
\int  ~ d\omega  dq     ~ 
 {\bf 1}_{|p|+|q| \le r(c)}
~ \mathcal{K}_{c} ~ [S_{c}(t) f(p)- S(t) f(p)]S_{c}(t)f_{c}(q)
\\
+
\int  ~ d\omega  dq     ~ 
 {\bf 1}_{|p|+|q| \le r(c)}
~ \mathcal{K}_{c} ~ [S_{c}(t) f_{c}(q)-S_{c}(t) f(q)] S(t) f(p)
\\
+
\int ~ d\omega  dq     ~ 
 {\bf 1}_{|p|+|q| \le r(c)}
~ \mathcal{K}_{c} ~ [S_{c}(t) f(q)-S(t) f(q)] S(t) f(p).
\end{multline*}
We will now estimate the first and second term in full detail, and we will then explain how the third term can be estimated in exactly the same way as the second term.

To estimate term $\mathbb{I}$, we use the following bound from Lemma \ref{asymptoticBB}
$$
\left| \mathcal{K}_{c}  -   \mathcal{K}_{\infty} \right|
\le A (1+|p|+|q|)^9/c^2.
$$
For simplicity without loss of generality we have fixed $m=n=9$. Thus we have
\begin{gather*}
\| \mathbb{I} \|_{L^1_p L^\infty_x}
\le 
Ab_2^2 r(c)^9/c^2
\int  ~ d\omega dq dp     ~ 
\mu^\beta(p) \mu^\beta(q) 
\le 
A b_2^2r(c)^9/c^2.
\end{gather*}
This follows from the uniform estimate $\| f^{\#}\|_{\infty}\le b_2$
assumed in Theorem \ref{NewtonianLimitThm}.

Next we estimate each of the terms in $\mathbb{II}$.
With 
the Glassey-Strauss variables \eqref{postCOLLvelGS},
 \eqref{omegaREP} and \eqref{postCOLLest} when ${|p|+|q| \le r(c)}$ we have
for some $A>0$ that
$$
|p^\prime | \le \left| p^{\prime} - \bar{p}^{\prime}  \right| +  \left| \bar{p}^{\prime}  \right|
\le A \frac{r(c)^3}{c^2} + 2\left( |p|+|q| \right) .
$$
A similar estimate also holds for $|q^\prime |$.  We conclude that
${|p|+|q| \le r(c)}$ implies
$$
|p'| + |q'| \le A r(c).
$$
We will use this argument implicitly several times in the following estimates.

For \eqref{t9}, we use a pre-post collisional change of variable  and end up with an upper bound in terms of our unknown difference.
After integration over $dp$, we have
\begin{multline*}
\left| 
\int  ~ d\omega  dq dp  ~ 
 {\bf 1}_{|p|+|q| \le r(c)}
~  \mathcal{K}_{c} ~ [S_{c}(t) f_{c}(p^{\prime})- S_{c}(t) f(p^{\prime})]S_{c}(t)f_{c}(q^{\prime})
\right|
\\
\le 
\int  ~ d\omega  dq dp  ~ 
 {\bf 1}_{|p'|+|q'| \le A r(c)}
~  \mathcal{K}_{c} ~ \|  f_{c}(p^{\prime})-  f(p^{\prime})\|_{L^\infty_x} \| f_{c}(q^{\prime})\|_{L^\infty_x}
\\
\le 
\int  ~ d\omega  dq dp  ~ 
 {\bf 1}_{|p|+|q| \le A r(c)}
~  \mathcal{K}_{c}(p,q,\omega) ~ \|  f_{c}(p)-  f(p)\|_{L^\infty_x} \| f_{c}(q)\|_{L^\infty_x}
\\
\le 
A \sigma_0  r(c) 
\| f_c - f \|_{L^1_p L^\infty_x} b_2.
\end{multline*}
As usual $b_2$ is the constant from Theorem \ref{NewtonianLimitThm}.  We have also used the estimate \eqref{kernelUPPERbdsC} combined with the exponential decay in $J^\beta(q)$ above.

We look now at the difference  \eqref{t9a}, which is a difference in the
the spatial argument.  We expand it out as follows
\begin{gather*}
S_{c}(t) f(p^{\prime})- S(t) f(p^{\prime})
=
f(t,x-t\hat{p},p^{\prime})
- 
f(t,x-t p,p^{\prime})
\\
=
f(t,x-t p +t(\hat{p}-p),p^{\prime})
- 
f(t,x-t p,p^{\prime})
\\
=
\tau_u^x S(t) f(p^{\prime})
- 
S(t) f(p^{\prime}),
\quad
u=t(\hat{p}-p).
\end{gather*}
On $|p|\le r(c)$, with \eqref{diffPhat}, $0\le t\le T$, and $c \ge 1$, we have
\begin{gather}
\label{uEST}
\left| u \right|
=
\left| t(\hat{p} - p) \right|
=
t \left| \frac{p}{\sqrt{1+|p|^2/c^2}} - p \right|
\le
\frac{T  }{c^2} (1+|p|^3)
\le
2\frac{T  }{c^2} r(c)^3.
\end{gather}
Taking the supremum in the $x$ variable, integrating over $dp$ and putting all of the estimates in this paragraph into  \eqref{t9a}, we obtain the upper bound for $\| \eqref{t9a} \|_{L^1_p L^\infty_x}$ 
\begin{multline*}
\le
\int  ~ d\omega  dq dp  ~ 
 {\bf 1}_{|p|+|q| \le r(c)}
~  \mathcal{K}_{c} ~ \|S_{c}(t) f(p^{\prime})- S(t) f(p^{\prime})\|_{L^\infty_x}
\|S_{c}(t)f_{c}(q^{\prime})\|_{L^\infty_x}
\\
\le
\sup_{|u|\le \frac{T  }{c^2} r(c)^3} \int  ~ d\omega  dq dp  ~ 
 {\bf 1}_{|p'|+|q'| \le A r(c)}
~  \mathcal{K}_{c} ~ \|\tau_u^x f(p^{\prime})- f(p^{\prime})\|_{L^\infty_x}
\|f_{c}(q^{\prime})\|_{L^\infty_x}
\\
\le
\sup_{|u|\le \frac{T  }{c^2} r(c)^3} \int  ~ d\omega  dq dp  ~ 
 {\bf 1}_{|p|+|q| \le A r(c)}
~  \mathcal{K}_{c} 
~ \|\tau_u^x f(p)- f(p)\|_{L^\infty_x}
\|f_{c}(q)\|_{L^\infty_x}.
\end{multline*}
We have done a pre-post collisional change of variable.
Furthemore
\begin{gather*}
\le 
A
\sigma_0 b_2 r(c)
\sup_{|u|\le \frac{T  }{c^2} r(c)^3}
\int  ~  dp  ~ 
 ~  \|\tau_u^x f(p)- f(p)\|_{L^\infty_x}
\\
\le 
A
\sigma_0 r(c)
\sup_{|u|\le \frac{T  }{c^2} r(c)^3}
\|\tau_u^x f- f\|_{L^1_p L^\infty_x}.
\end{gather*}
These last estimates follow from the upper bound for $f_c$ in Theorem \ref{NewtonianLimitThm}, combined with the estimate for $\mathcal{K}_c$ in \eqref{kernelUPPERbdsC}, and also Lemma \ref{JSOLbounds}.

For \eqref{t10a}, we use the method from \eqref{t9}. 
Then for \eqref{t11}, we use the method from \eqref{t9a}, except now we use the symmetric reverse estimate for the kernel $\mathcal{K}_{c}(p,q,\omega)$ as given in the display just below \eqref{kernelUPPERbdsC}. 

For \eqref{t12}, we estimate 
$\left| \bar{q}^{\prime} - q^{\prime} \right|$ as in
\eqref{postCOLLest}.
 With that, we have
\begin{gather*}
S(t) f(q^{\prime})- S(t) f(\bar{q}^{\prime})
=
S(t) f(q^{\prime})- S(t) f(q^{\prime}+(\bar{q}^{\prime} - q^{\prime}))
\\
=
S(t) f(q^{\prime})- \tau^p_{u} S(t) f(q^{\prime}),
\quad
u=\bar{q}^{\prime} - q^{\prime}.
\end{gather*}
Now the rest of the estimate for \eqref{t12} can be handled exactly as in the estimate for \eqref{t9a}, the only difference being that the difference is in terms of $\tau^p_{u}$ rather than $\tau^x_{u}$.  The end result is that 
\begin{gather*}
\| \eqref{t12} \|_{L^1_p L^\infty_x}
\le 
A
\sigma_0 r(c)
\sup_{|u|\le \frac{T  }{c^2} r(c)^3}
\|\tau_u^p f- f\|_{L^1_p L^\infty_x}.
\end{gather*}
The remaining term \eqref{t10} needs to be estimated in a slightly different way.

Lastly we estimate the term \eqref{t10} which forms
a difference in the momentum arguments of the same function.  For this term we expand the differences as 
\begin{gather*}
S(t) f(p^{\prime})- S(t) f(\bar{p}^{\prime})
=
S(t) f(\bar{p}^{\prime}+(p^{\prime} - \bar{p}^{\prime}))- S(t) f(\bar{p}^{\prime})
\\
=
\tau^p_{\bar{u}}S(t) f(\bar{p}^{\prime})-  S(t) f(\bar{p}^{\prime}),
\quad
\bar{u}=  p^{\prime} - \bar{p}^{\prime}.
\end{gather*}
On $|p|+|q|\le r(c)$, with  \eqref{postCOLLest}, and $c \ge 1$, we have
$$
\left| \bar{u} \right|
=
\left| p^{\prime} - \bar{p}^{\prime} \right|
\le
A \frac{r(c)^3  }{c^2} .
$$
Taking the supremum in the $x$ variable, integrating over $dp$ and putting all of the estimates in this paragraph into  \eqref{t10}, we obtain the upper bound
\begin{gather*}
\| \eqref{t10} \|_{L^1_p L^\infty_x}
\le
\int  ~ d\omega  dq dp  ~ 
 {\bf 1}_{|p|+|q| \le r(c)}
~  \mathcal{K}_{c} ~ \|S(t) f(p^{\prime})- S(t) f(\bar{p}^{\prime})\|_{L^\infty_x}
\|S(t)f(\bar{q}^{\prime})\|_{L^\infty_x}
\\
\le
\sup_{|\bar{u}|\le A' \frac{r(c)^3  }{c^2} } \int  ~ d\omega  dq dp  ~ 
 {\bf 1}_{|p|+|q| \le r(c)}
~  \mathcal{K}_{c}  ~ \|\tau_{\bar{u}}^p f(\bar{p}^{\prime})- f(\bar{p}^{\prime})\|_{L^\infty_x}
\|f(\bar{q}^{\prime})\|_{L^\infty_x}
\\
\le
\sigma_0 r(c)^{7/2} \sup_{|\bar{u}|\le A' \frac{r(c)^3  }{c^2} } \int  ~ d\omega  dq dp  ~ 
 {\bf 1}_{|p|+|q| \le r(c)}
 ~ \|\tau_{\bar{u}}^p f(\bar{p}^{\prime})- f(\bar{p}^{\prime})\|_{L^\infty_x}
\|f(\bar{q}^{\prime})\|_{L^\infty_x}
\\
\le
\sigma_0 r(c)^{7/2}
\sup_{|\bar{u}|\le A' \frac{r(c)^3  }{c^2} } \int  ~ d\omega  dq dp  ~ 
 {\bf 1}_{|p|+|q| \le 4 r(c)}
~ \|\tau_{\bar{u}}^p f(p)- f(p)\|_{L^\infty_x} ~ 
\mu^{\beta/2}(q).
\end{gather*}
We have used \eqref{kernelUPPERbdsC},
and the pre-post change  in the Newtonian variables \eqref{omegaREP}.
Then
\begin{gather*}
\| \eqref{t10} \|_{L^1_p L^\infty_x}
\le 
A
\sigma_0 r(c)^{7/2}
\sup_{|\bar{u}|\le A' \frac{r(c)^4  }{c^2} }
\|\tau_{\bar{u}}^p f- f\|_{L^1_p L^\infty_x}.
\end{gather*}
This completes all of our estimates for the terms in $\mathbb{II}$.

Notice that each of the terms in $\mathbb{III}$ has an analogous term in $\mathbb{II}$, except there are not post-collisional velocities to be concerned with.  Thus the estimates in $\mathbb{III}$ follow exactly as the estimates in 
$\mathbb{II}$ with the same result.

We combine all of the estimates for each of the terms $\mathbb{I}$,
$\mathbb{II}$,
and
$\mathbb{III}$
to obtain
\begin{multline*}
\| S_{c}(t)\mathcal{Q}_{c}(t) -  S(t)\mathcal{Q}(t)\|_{L^1_p L^\infty_x}
\le
A_1 r(c) 
\| f_c - f \|_{L^1_p L^\infty_x} 
\\
+
A_2 r(c)^{7/2}
\sup_{|u|\le A_3 r(c)^3/c^2  } \left(\|\tau_u^x f- f\|_{L^1_p L^\infty_x}
+
\|\tau_u^p f- f\|_{L^1_p L^\infty_x}
\right)
+
A_4 \frac{r(c)^9}{c^{2}}.
\end{multline*}
From this inequality and 
Lemma  \ref{controlDIFF} 
we 
conclude the statement of Lemma \ref{nonLINasy}
 for any small $\varepsilon >0$,
since $r(c) = \left(\log c\right)^{\alpha_1}$ with $1/2 < \alpha_1 <1$.
\end{proof}

Lastly, we  prove the main Newtonian Limit Theorem \ref{NewtonianLimitThm}.  

\begin{proof}
We consider the difference of solutions to 
\eqref{mildSOL} 
and \eqref{mildSOLnewt}
with $ c \ge 1$ as:
\begin{multline*}
f_{c}(t,x,p) - f(t,x,p) = S_{c}(t)f_{0,c}(x,p) - S(t)f_{0}(x,p) 
\\
+ \int_0^t ~ ds ~ S_{c}(t-s)Q_{c}(s, x, p) - S(t-s)Q(s, x, p).
\end{multline*}
After taking the $L^1_p L^\infty_x$ norm we obtain the following inequality
\begin{multline*}
\| f_{c}(t) - f(t) \|_{L^1_p L^\infty_x}  
\le  
\|S_{c}(t)f_{0,c} - S(t)f_{0}  \|_{L^1_p L^\infty_x}  
\\
+ 
\int_0^t ~ ds ~\| S_{c}(t-s)\mathcal{Q}_{c}(s) -  S(t-s)\mathcal{Q}(s)\|_{L^1_p L^\infty_x}.
\end{multline*}
By translation invariance of the $L^\infty_x$ norm, we have
\begin{multline*}
\|S_{c}(t)f_{0,c} - S(t)f_{0} \|_{L^1_p L^\infty_x}  
\le
\| S_{c}(t)f_{0,c} - S_{c}(t)f_{0} \|_{L^1_p L^\infty_x}
+
\| S_{c}(t)f_{0} - S(t)f_{0} \|_{L^1_p L^\infty_x}
\\
=
\| f_{0,c} - f_{0} \|_{L^1_p L^\infty_x}
+
\|\tau^x_{u} S(t)f_{0} - S(t)f_{0} \|_{L^1_p L^\infty_x},  
\end{multline*}
with ${u} \eqdef  t(\hat{p} - p)$.  
When $|p| \le r(c)$,
and
$r(c) = \left(\log c\right)^{\alpha_1}$ with $1/2 < \alpha_1 <1$, 
we use the estimate in \eqref{uEST}
to conclude
$$
|u|
\le 
AT/c^{2-\varepsilon},
\quad
\forall \varepsilon\in(0,2].
$$
Thus under the assumptions of Theorem \ref{NewtonianLimitThm}, on $|p| \le r(c)$,   we have
\begin{gather*}
\sup_{|u| \le AT/c^{2-\epsilon}}
\|\tau^x_{u} S(t)f_{0} - S(t)f_{0} \|_{L^1_p L^\infty_x}  
\le
A(T)/c^{2-\varepsilon},
\quad
\forall \varepsilon\in(0,2].
\end{gather*}
On the other hand, if $|p| \ge r(c)$, then we use estimates in the spirit of \eqref{extraCest}:
\begin{multline}
\notag
\| S_{c}(t)f_{0} - S(t)f_{0} \|_{L^1_p(|p| \ge r(c)) L^\infty_x}
\le
2\| f_{0} \|_{L^1_p(|p| \ge r(c)) L^\infty_x}
\\
=
2
\int_{\mathbb{R}^N} ~dp ~ {\bf 1}_{|p|\ge r(c)} ~ \|f_{0}(p)\|_{L^\infty_x}
\le
2b_2
\int_{\mathbb{R}^N} ~dp ~ {\bf 1}_{|p|\ge r(c)} ~ \mu^\beta (p)
\\
\le 
A ~ \mu^{\beta/2}(r(c)) 
\int_{\mathbb{R}^N} ~ dp ~ {\bf 1}_{|p| > r(c)}   ~ \mu^{\beta/2}(p)
\\
\le 
A  \mu^{\beta/2}(r(c)) 
\le
 A ~ e^{-\beta (r(c))^2/4}
 \le 
 \frac{A}{c^{2}}. 
\end{multline}
The estimates above follow from the bounds on the limit solution in Theorem \ref{NewtonianLimitThm} and \eqref{weightINF}.
Finally the non-linear asymptotic estimate from Lemma \ref{nonLINasy} yields
\begin{multline*}
\| S_{c}(t-s)\mathcal{Q}_{c}(s) -  S(t-s)\mathcal{Q}(s)\|_{L^1_p L^\infty_x}
\\
\le
A ~ r(c) ~ 
\| f_{c}(s) - f(s) \|_{L^1_p L^\infty_x} 
+
A(T) / c^{2-\varepsilon},
\quad 
\forall \varepsilon\in(0,2].
\end{multline*}
Putting all of these estimates together we have
\begin{multline*}
\| f_{c}(t) - f(t) \|_{L^1_p L^\infty_x}  
\le  
\| f_{0,c} - f_{0} \|_{L^1_p L^\infty_x}+
A(T)/c^{2-\varepsilon}
\\
+ 
A ~ r(c) ~ 
\int_0^t ~ ds ~
\| f_{c}(s) - f(s) \|_{L^1_p L^\infty_x}.  
\end{multline*}
By the Gronwall inequality
$$
\| f_{c}(t) - f(t) \|_{L^1_p L^\infty_x}  
\le
\left(
 \| f_{0,c} - f_{0} \|_{L^1_p L^\infty_x}+
A(T)/c^{2-\varepsilon}
\right)
e^{A T  r(c)}.
$$
By the slow growth of $r(c)$, with \eqref{slowC}, this gives easily the Cauchy property with the convergence rate of $1/c^{k-\delta}$ as stated in Theorem \ref{NewtonianLimitThm}.
\end{proof}

This last proof completes the main focus of this paper.   We conclude with two Appendicies.  In Appendix A, we elucidate a collection of Lorentz Transformations that may be useful for choosing coordinates as in \eqref{postC}. 
In Appendix B, we point out some differential cross-sections which appear in the physics literature for the relativistic Boltzmann equation.

\section*{Appendix A:  Lorentz Transformations}

In this appendix, we will write down three Lorentz Transformations which may be useful in choosing coordinates in relativistic Kinetic Theory.  
In the prior sections, we worked in arbitrary dimensions $N \ge 2$.  For simplicity, in this section we work in the physical dimension $N =3$.
In this section we will also use the notation
$$
P=\left(\begin{array}{c} p^{0} \\ p\end{array}\right), 
\quad 
Q=\left(\begin{array}{c} q^{0} \\ q\end{array}\right). 
$$
 For more on Lorentz Transformations, we refer to \cite{weinbergBK,naberBK}.  
 We give their definition:

\begin{definition}\label{rcop:LTdef} A four by four matrix
$\Lambda$ is a  Lorentz Transformation if
\begin{equation*}
(\Lambda P)\cdot (\Lambda Q)=P\cdot Q.
\end{equation*}
Let $D=\text{diag}(1, -1, -1, -1)$.  Equivalently, $\Lambda$ is a Lorentz Transformation if 
\begin{equation*}
\Lambda^T D\Lambda =D.
\end{equation*}
\end{definition}

Any Lorentz Transformation, $\Lambda$, is invertible and the last display implies 
$$
\Lambda^{-1}=D \Lambda^T D.
$$
Thus the inverse given in \eqref{postC} can be easily computed.  Since $\Lambda$ has sixteen components and is restricted by ten equations,  $\Lambda$  has six free parameters.

Perhaps the most well known Lorentz transformation is the boost matrix.  Given 
$
v=(v^1,v^2,v^3)\in\mathbb{R}^3,
$  
we write the boost matrix as
$$
\Lambda_{b}
=
\left(
\begin{array}{cccc}
\gamma & -\gamma v^1 & -\gamma v^2 & -\gamma v^3
\\
-\gamma v^1 & 1+(\gamma-1)\frac{v^1 v^1}{|v|^2} & (\gamma-1)\frac{v^2 v^1}{|v|^2} & (\gamma-1)\frac{v^3 v^1}{|v|^2}
\\
-\gamma v^{2} & (\gamma-1)\frac{v^1 v^2}{|v|^2} & 1+(\gamma-1)\frac{v^2 v^2}{|v|^2} & (\gamma-1)\frac{v^3 v^2}{|v|^2}
\\
-\gamma v^{3} & (\gamma-1)\frac{v^1 v^3}{|v|^2} & (\gamma-1)\frac{v^2 v^3}{|v|^2} & 1+(\gamma-1)\frac{v^3 v^3}{|v|^2}
\end{array}
\right),
$$
where $\gamma=(1-|v|^2)^{-1/2}$.  Notice that $\Lambda_{b}$ has only three free parameters.  It is well-known that any Lorentz Transformation can be expressed as a rotation multiplied by a boost matrix.

We are exclusively concerned with proper orthochronous Lorentz transformations $\Lambda$ which send $P+Q$ into the center-of-momentum system
\begin{equation}
\Lambda(P+Q)=(\sqrt{s}, 0, 0, 0)^t.
\label{lw:lorentzS}
\end{equation}
We recall that $s$ is defined in \eqref{sDEFINITION}.
We will give three examples of Lorentz Transformations satisfying \eqref{lw:lorentzS}.  \\

\noindent{\it Example 1: The Boost Matrix}.   We choose $v$ such that $\Lambda_{b}$ satisfies (\ref{lw:lorentzS}).  Let
$$
v=\frac{p+q}{p_0+q_0}, \quad \gamma=\frac{p_0+q_0}{\sqrt{s}}.
$$
Then  $\Lambda_{b}$  is given by
$$
\Lambda_{b}
=
\left(
\begin{array}{cccc}
\frac{p_0+q_0}{\sqrt{s}} & -\frac{p_1+q_1}{\sqrt{s}} & -\frac{p_2+q_2}{\sqrt{s}} & -\frac{p_3+q_3}{\sqrt{s}}
\\
-\frac{p_1+q_1}{\sqrt{s}} 
& 1+\left(\gamma-1\right)\frac{v^1 v^1}{|v|^2} 
&  \left(\gamma-1\right)\frac{v^2 v^1}{|v|^2}
& \left(\gamma-1\right)\frac{v^3 v^1}{|v|^2}
\\
-\frac{p_2+q_2}{\sqrt{s}} 
& \left(\gamma-1\right)\frac{v^1 v^2}{|v|^2} 
& 1+\left(\gamma-1\right)\frac{v^2 v^2}{|v|^2} 
& \left(\gamma-1\right)\frac{v^3 v^2}{|v|^2}
\\
-\frac{p_3+q_3}{\sqrt{s}} 
& \left(\gamma-1\right)\frac{v^1 v^3}{|v|^2} 
& \left(\gamma-1\right)\frac{v^2 v^3}{|v|^2} 
& 1+\left(\gamma-1\right)\frac{v^3 v^3}{|v|^2}
\end{array}
\right),
$$
where
$
\frac{v^i v^j}{|v|^2} =\frac{(p_i+q_i)(p_j+q_j)}{|p+q|^2}.
$
By a direct calculation, this example satisfies (\ref{lw:lorentzS}).  We notice that $\Lambda_b \to I_4$, the four by four identity matrix, as $c\to\infty$. \\

 \noindent{\it Example 2}.  
Given $p,q\in \mathbb{R}^N$ with $p+q\ne 0$, we can always choose three orthonormal vectors
$$
w^1, \quad w^2, \quad w^3=\frac{p+q}{|p+q|}.
$$
If, for example, $p_1+q_1\ne 0$ and $p_2+q_2\ne 0$, then we can explicitly write
\begin{eqnarray*}
w^1&=&\frac{\left(-(p_2+q_2)(p_3+q_3), 2(p_1+q_1)(p_3+q_3), -(p_2+q_2)(p_3+q_3)\right)}
{|\left(-(p_2+q_2)(p_3+q_3), 2(p_1+q_1)(p_3+q_3), -(p_2+q_2)(p_3+q_3)\right)|}
\\
w^2&=&w^1\times w^3/|w^1\times w^3|.
\end{eqnarray*}
If instead $p_1+q_1=0$, then
\begin{eqnarray*}
w^1=(1,0,0), \quad
w^2=w^1\times w^3/|w^1\times w^3|.
\end{eqnarray*}
With this notation, we can write down a Lorentz Transformation which maps $p+q\to 0$ as
\begin{equation}
\Lambda=
\left(
\begin{array}{cccc}
\frac{p_0+q_0}{|P+Q|}  & -\frac{p_1+q_1}{|P+Q|} & -\frac{p_2+q_2}{|P+Q|}  & -\frac{p_3+q_3}{|P+Q|} 
\\
0 & w^1_{1} & w^1_{2} & w^1_{3}
\\
0 & w^2_{1} & w^2_{2} & w^2_{3}
\\
\frac{|p+q|}{|P+Q|}& -\frac{p_1+q_1}{|p+q|}\frac{p_0+q_0}{|P+Q|} & 
-\frac{p_2+q_2}{|p+q|}\frac{p_0+q_0}{|P+Q|}  & -\frac{p_3+q_3}{|p+q|}\frac{p_0+q_0}{|P+Q|} 
\end{array}
\right).
\label{lw:LT}
\end{equation}
This Lorentz transformation clearly satisfies (\ref{lw:lorentzS}).  This example is fully derived in \cite{strainPHD}. \\

\noindent{\it Example 3: Hilbert-Schmidt form for the Linearized Operator}.\label{rcop:HSlorentz}
We now wish to derive a Lorentz Transformation which satisfies \eqref{lw:lorentzS} but also
 \begin{equation}
\Lambda(P-Q)=(0, 0, 0, g)^t.
\label{lw:lorentzG}
\end{equation}
In \cite[p.277]{MR635279}  this transformation is used to write down a Hilbert-Schmidt form for the linearized collision operator.  But the transformation is not written down explicitly therein.  What is given is conditions \eqref{lw:lorentzS}, \eqref{lw:lorentzG} and the first row $\Lambda_0$.  From this information we can make the following educated guess for $\Lambda$ satisfying (\ref{lw:lorentzS}) and (\ref{lw:lorentzG}):
$$
\Lambda
=
\left(
\begin{array}{cccc}
\frac{p_0+q_0}{\sqrt{s}} & -\frac{p_1+q_1}{\sqrt{s}} & -\frac{p_2+q_2}{\sqrt{s}} &-\frac{p_3+q_3}{\sqrt{s}}
\\
\Lambda_{0}^{1} & \Lambda_{1}^{1} & \Lambda_{2}^{1} & \Lambda_{3}^{1}
\\
0 & \frac{(p\times q)_1}{|p\times q|}& \frac{(p\times q)_2}{|p\times q|}& \frac{(p\times q)_3}{|p\times q|}
\\
\frac{p_0-q_0}{g} & -\frac{p_1-q_1}{g} & -\frac{p_2-q_2}{g} &-\frac{p_3-q_3}{g}
\end{array}
\right).
$$
 By the Lorentz condition in Definition \ref{rcop:LTdef}, we may  conclude 
$$
\Lambda_{0}^1= \frac{2|p\times q|}{\sqrt{s}g}=\frac{|p\times q|}{\sqrt{(P\cdot Q)^2-c^4}}.
$$
The coefficients $\Lambda^{1}_{i}$ ($i=1,2,3$) are similarly completely determined by the Lorentz condition.  However, this turns out to be a long computation.  In the end we obtain
$$
\Lambda^{1}_{i}=
\frac{2\left( p_i\left\{p_0-q_0P\cdot Q\right\}+q_i\left\{q_0-p_0P\cdot Q\right\}\right)}{\sqrt{s}g|p\times q|}, 
\quad (i=1,2,3).
$$
Now, we have a complete description of this Lorentz transformation in terms of $p, q$.  

The way we derived this matrix was by guessing the first, third and fourth rows and then using the Lorentz condition.  However, it is much easier to verify that this is indeed a Lorentz Transformation after the fact by checking  the condition in Definition \ref{rcop:LTdef}.  The interested reader may find this calculation in \cite{strainPHD}.

\section*{Appendix B:  Examples of relativistic Boltzmann cross-sections}  

In \cite{MR0156656}, Grad gave several examples of important collision kernels to study in the classical theory.  The same type of discussion would be important for the relativistic kinetic theory, in this direction we point out \cite{DEnotMSI}.  Here we  write down 
 a few scattering kernels which can be found in the physics literature.

The calculation of the differential cross-section in the relativistic situation utilzes quantum field theory, see for instance \cite{MR1402248}.  Cross sections are not presently derived from a scattering problem (as it is done in the Newtonian case) in particular because there is no widely accepted theory of relativistic N-Body dynamics at present.

Above, we have written everything down with the mass $m$ normalized to one, $m=1$.  Here we write down the mass without normalization.

\subsubsection*{Short Range Interactions}  \cite{MR1402248,MR1958975}.  For short range interactions, 
\begin{equation}
\sigma\eqdef \text{constant}.   
\label{hardSPHERE}
\end{equation}
This ``hard-ball'' cross section is the relativistic analouge of the hard-sphere kernel in the Newtonian case.  Indeed, as we have already mentioned, the Newtonian limit of the relativistic Boltzmann equation in this case is the hard-sphere Boltzmann equation.   This is also sometimes written as $s\sigma = \text{constant}$.

\subsubsection*{M{\o}ller Scattering}   \cite[p.350]{MR635279}. M{\o}ller scattering is used as an approximation of electron-electron scatttering.
$$
\sigma(g,\theta)=r_0^2\frac{1}{u^2(u^2-1)^2}\left\{\frac{(2u^2-1)^2}{\sin^4\theta}
-
\frac{2u^4-u^2-\frac{1}{4}}{\sin^2\theta}
+\frac{1}{4}(u^2-1)^2 \right\},  
$$
where the magnitude of total four-momentum scaled w.r.t. total mass is
$$
u=\frac{\sqrt{s}}{2mc},
$$
and $r_0=\frac{e^2}{4\pi mc^2}$ is the classical electron radius.  Photon-photon scattering is often neglected because the size of the cross-section is `negligible'.

\subsubsection*{Compton Scattering}      \cite[p.351]{MR635279}.   Compton scattering is an approximation of photon-electron scattering.  
$$
\sigma(g,\theta)=\frac{1}{2}r_0^2(1-\xi)\left\{1
+
\frac{1}{4}\frac{\xi^2 (1-\cos\theta)^2}{1-\frac{1}{2}\xi(1-\cos\theta)}
+
\left(
\frac{1-(1-\frac{1}{2}\xi)(1-\cos\theta)}{1-\frac{1}{2}\xi(1-\cos\theta)}
\right)^2\right\}, 
$$
where
$$
\xi=1-\frac{m^2c^2}{s}.
$$
See \cite[p.81]{MR1958975} for the Newtonian limit in this case.

\subsubsection*{(elastic) Neutrino Gas}  \cite[p.478]{MR0471665}  For a neutrino gas, the differential cross section is independent of the scattering a angle $\theta$.   
$$
\sigma(g,\theta)=\frac{G^2}{\pi \hbar^2 c^2} g^2,
$$
where $G$ is the weak coupling constant and $2\pi \hbar$ is Plank's constant.  These are massless particles.  See also \cite[p.290]{MR635279}.

\subsubsection*{Israel particles}  \cite[p.1173]{MR0165921}    The Isreal particles are the analogue of the well known ``maxwell molecules'' cross section in the Newtonian theory.  They are 
$$
\sigma(g,\theta)= \frac{m}{2g}\frac{ b(\theta)}{1+\left(g/mc\right)^2}.
$$
With this cross section, Israel derived eigenfunctions for the Linearized relativistic Boltzmann collision operator.  Variants  are used in \cite{polakRKGT,MR1898707}.  For instance the cross section for ``Maxwell Particles" is formed by removing the factor $1+\left(g/mc\right)^2$.  It converges to the maxwell molecules cross section in the Newtonian limit.  


 \subsection*{Acknowledgments}
The author would like to thank  Simone Calogero, Robert Glassey, Seung-Yeal Ha, David Levermore, and Nader Masmoudi for helpful comments and discussions after a recent talk at IPAM on this paper.  
  The author also thanks the anonymous referees for their helpful comments.

\end{document}